\documentclass[preprint,12pt,authoryear]{elsarticle}

\usepackage{amssymb}

\usepackage[utf8]{inputenc} 
\usepackage[T1]{fontenc}    
\usepackage[colorlinks=true,linkcolor=blue,citecolor=blue,urlcolor=blue]{hyperref}
\usepackage{url}            
\usepackage{booktabs}       
\usepackage{amsfonts}       
\usepackage{nicefrac}       
\usepackage{microtype}      
\usepackage{lipsum}
\usepackage{fancyhdr}       
\usepackage{graphicx}       
\usepackage{amsmath}
\usepackage{mathtools}
\usepackage{algorithm}
\usepackage{algpseudocode}
\graphicspath{{media/}}     
\usepackage{adjustbox}
\usepackage{bibentry}
\usepackage{verbatim}
\usepackage{geometry}

\usepackage{longtable}
\usepackage[colorinlistoftodos,textwidth=1in]{todonotes}
\usepackage{placeins}

\newcommand{\bigM}{\mathbf{M}}
\newcommand{\bigQ}{\mathbf{Q}}
\newcommand{\bigt}{\mathbf{t}}

\newcommand{\bigy}{\mathbf{y}}
\newcommand{\bigphi}{\boldsymbol{\phi}}
\newcommand{\bigtheta}{\boldsymbol{\theta}}
\newcommand{\R}{\mathbb{R}}
\newcommand{\bigq}{\mathbf{q}}
\newcommand{\bigSigma}{\mathbf{\Sigma}}
\newcommand{\trace}{\text{tr}}
\newcommand{\Var}{\text{Var}}
\newcommand{\Cov}{\text{Cov}}

\journal{Computers and Chemical Engineering}

\begin{document}

\begin{frontmatter}



\title{Measure This, Not That: Optimizing the Cost and Model-Based Information Content of Measurements}


\author[inst1]{Jialu Wang}

\affiliation[inst1]{organization={Department of Chemical and Biomolecular Engineering},
            addressline={University of Notre Dame}, 
            city={Notre Dame},
            state={Indiana 46556},
            country={United States}}

\author[inst2]{Zedong Peng}
\author[inst3,inst4]{Ryan Hughes}
\author[inst5]{Debangsu Bhattacharyya}
\author[inst2]{David E.~Bernal Neira}
\author[inst1]{Alexander W.~Dowling}

\affiliation[inst2]{organization={Davidson School of Chemical Engineering},
            addressline={Purdue University}, 
            city={West Lafayette},
            state={Indiana 47907},
            country={United States}}

\affiliation[inst3]{organization={National Energy Technology Laboratory},
            city={Pittsburgh},
            state={Pennsylvania 15236},
            country={United States}}

\affiliation[inst4]{organization={Support Contractor, National Energy Technology Laboratory},
            city={Pittsburgh},
            state={Pennsylvania 15236},
            country={United States}}

\affiliation[inst5]{organization={Department of Chemical and Biomedical Engineering},
            addressline={West Virginia University}, 
            city={Morgantown},
            state={West Virginia 26506},
            country={United States}}

\begin{abstract}

Model-based design of experiments (MBDoE) is a powerful framework for selecting and calibrating science-based mathematical models from data. This work extends popular MBDoE workflows by proposing a convex mixed integer (non)linear programming (MINLP) problem to optimize the selection of measurements. The solver \texttt{MindtPy} is modified to support calculating the D-optimality objective and its gradient via an external package, \texttt{SciPy}, using the grey-box module in \texttt{Pyomo}. The new approach is demonstrated in two case studies: estimating highly correlated kinetics from a batch reactor and estimating transport parameters in a large-scale rotary packed bed for CO$_2$ capture. Both case studies show how examining the Pareto-optimal trade-offs between information content measured by A- and D-optimality versus measurement budget offers practical guidance for selecting measurements for scientific experiments.
\end{abstract}

\begin{graphicalabstract}
\includegraphics{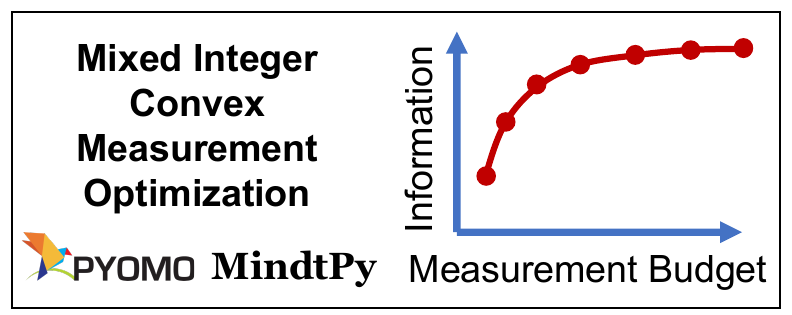}
\end{graphicalabstract}

\begin{highlights}
\item Proposes convex MI(N)LP problems to optimize the Fisher information of measurements
\item Utilizes grey-box module in \texttt{Pyomo} to incorporate externally evaluated objective function
\item Optimizes with highly correlated measurements on kinetics case study
\item Demonstrates scalability in CO$_2$ capture case study
\end{highlights}


\end{frontmatter}


\section{Introduction}

Design of experiments (DoE) plays an important role in creating mathematical models by identifying the most informative data while minimizing the consumption of materials, time, and human resources \citep{franceschini2008model}. The transformative potential digital twins, i.e., computational models that mimic complex physical systems and are adaptively updated with new data \citep{rasheed2019digital, zobel2021digital, ors2020conceptual}, emphasize the importance of optimizing data collection when building and validating predictive models. Traditional black-box DoE methods, such as factorial design, response surface methods, and space-filling designs, are widely employed \citep{box1992experimental, myers2016response}, as they do not require a science-based mathematical model. In contrast, model-based design of experiments (MBDoE) uses science-based models, such as:
\vspace{-5pt}
\begin{equation}
    \bigy(t) = \mathbf{f}(\mathbf{z}(t), \bigtheta, \bigphi(t)) + \boldsymbol{\epsilon},
\label{eq:v1}
\end{equation}

\noindent where $\bigy$ is the vector of measurements that are corrupted with random error $\boldsymbol{\epsilon}$, $\mathbf{z}$ is the vector of state variables, $\bigtheta$ is the vector of (physically interpretable) uncertain model parameters, and $\mathbf{u}(t)$ is the vector of control variables. Experimental decisions $\bigphi(t) = (\mathbf{u}(t), \mathbf{z}(t_0), \mathbf{t})$ include the control variables, experimental initial conditions, and the choice of measurement time points. MBDoE determines the experiment decisions $\bigphi(t)$ that produce the measurements $\bigy$(t) that provide the most information about $\bigtheta$. These models are typically ordinary differential equations (ODEs) or partial differential equations (PDEs) with time-varying decisions $\bigphi(t)$, which makes selecting the best $\bigphi(t)$ an infinite-dimensional optimization problem. In this case, MBDoE campaigns are calculated using methods from optimal control, wherein $\bigphi(t)$ is discretized, and gradient-based optimization computes the best piecewise constant, linear, \emph{et cetera}, policy \citep{franceschini2008model}. Using a discretization with many time points results in a high-dimensional optimization problem. Black-box DoE techniques suffer from the curse of dimensionality as the data requirements grow rapidly with the number of input dimensions. Thus black-box DoE methods are especially well suited for steady-state system \citep{soepyan2018sequential} which are often low dimensional. MBDoE for dynamical systems overcomes these issues by exploiting the model structure with derivative-based nonlinear optimization algorithms. For example, \citet{waldron2020model} showed that MBDoE reduced the parameter uncertainty by up to 40\% compared to the factorial design. \citet{wang2022pyomo} used MBDoE to predict \emph{a priori} the information gained from modifying experimental procedures or adding new measurements. Recent applications of MBDoE include reaction kinetics \citep{knoll2022autonomous,wang2022pyomo, tillmann2023development, cenci2023exploratory}, thermodynamic models \citep{befort2023data}, additive manufacturing \citep{wang2023physics, shahmohammadi2020sequential}, membrane separation \citep{liu2022membrane}, crystallization process \citep{yuan2023combined}, electrochemical models for fuel cells \citep{kravos2021methodology}, and bioreactors \citep{liang2020model, kim2019efficient, saccardo2023model}.

\begin{figure}[!htb]
\centering
\includegraphics[width=\textwidth]{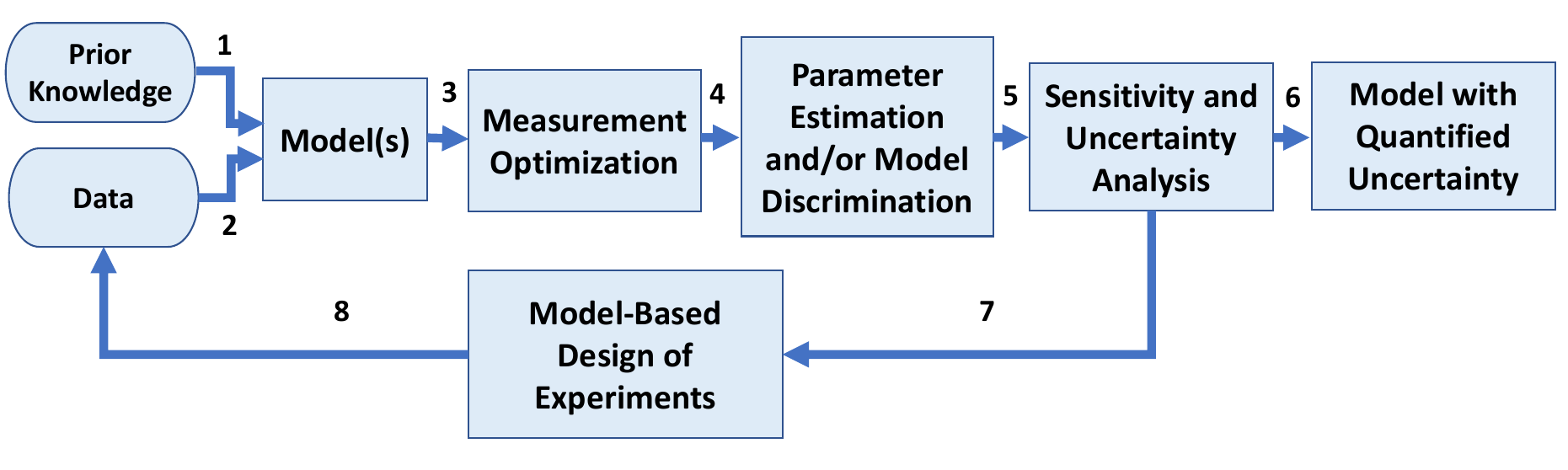}
\caption{Sequential MBDoE workflow extended from \citet{franceschini2008model}, \citet{ wang2022pyomo}, and references therein to consider measurement optimization. The workflow and numbered arrows are described in the text.} \label{fig:flowchart}
\end{figure}

A key challenge often overlooked in MBDoE literature is the economical selection of measurements.  
Fig.~\ref{fig:flowchart} shows the typical MBDoE workflow \citep{franceschini2008model, wang2022pyomo} extended to consider measurement optimization (MO). MBDoE starts by combining prior knowledge (arrow 1) with preliminary data (arrow 2) to postulate one or more mathematical models (arrow 3). These models are used to optimize measurement selection (arrow 4). Next, model discrimination and parameter estimation are performed to downselect the models and estimate parameters (arrow 5). Sensitivity and uncertainty analysis provides insights into the model (e.g., a parameter is not practically identifiable). If the parameter uncertainty is sufficiently small, the workflow terminates and outputs a model with quantified uncertainty (arrow 6). Otherwise, the models and parameter estimates (arrow 7) are used for MBDoE to recommend the next set of experiments (arrow 8), and the algorithm continues unless the experimental budget is exhausted.

In this workflow (Fig.~\ref{fig:flowchart}), MO decides which sensors to install. Different choices of measurements lead to different information contents, influencing the practical identifiability of the model parameters and deciding the significance of parameter estimates. However, there are many candidate measurement variables, e.g., sensor choices, locations, and sampling times in an experimental system. They can vary widely in terms of error variances and covariances, sampling times, sampling positions, and measurement techniques, resulting in a plethora of choices and constraints; it is often neither economically nor physically feasible to measure all state variables at all time points. 

\subsection{Optimization Strategies for MBDoE}

The MBDoE optimization problem is defined as: 
\vspace{-5pt}
\begin{equation}
    \bigphi^*(t) = \underset{\bigphi(t) \in \boldsymbol{\Phi}}{\text{argmax}} ~ \Psi(\bigM(\bigtheta, \bigphi(t)))
\label{eq:simple-opt}
\end{equation}

\noindent where $\bigphi(t)$ is the set of experimental design variables, which can be time-varying and therefore create a high-dimensional design space; $\boldsymbol{\Phi}$ denotes this feasible experiment design space, which often includes bounds on $\bigphi(t)$ and dynamic constraints on the states of the system; $\Psi(\cdot)$ is a design criterion; and $\bigM$ is the Fisher information matrix (FIM), which measures the experimental information content for estimating unknown parameters $\bigtheta$. 
Alternate objectives include model discrimination \citep{hunter1965designs, buzzi1983new, galvanin2016joint, tillmann2023development} and decreasing the uncertainties in model predictions \citep{kiefer1959optimum, dasgupta2021g, cenci2023exploratory}.

Standard MBDoE assumes the correct model structure, Gaussian priors for $\bigtheta$, Gaussian measurement error structure, and neglects parameter uncertainty. Several recent advances in MBDoE methods relax these assumptions. For example, robust MBDoE extensions \citep{telen2013guaranteed, mesbah2015probabilistic,petsagkourakis2021safe} consider uncertainty by optimizing the expected value, variance, or similar risk metric of the MBDoE objective(s). Online approaches \citep{galvanin2012online, kim2019efficient, pankajakshan2019multi, waldron2019closed} sequentially incorporate newly achieved knowledge from the newest experiments, update models, and solve the next MBDoE problem automatically so that experiments and MBDoE can be conducted autonomously. We emphasize the importance of automated MBDoE for the emerging fields of digital twins and self-driving laboratories \citep{lair2024critical,agi2024computational}. \citet{rodrigues2022tractable, lainez2015stochastic, shen2023bayesian} developed MBDoE with more sophisticated error distributions. \citet{reichert2019influence} studied the influence of various error descriptions in the measurement data for MBDoE approaches, including random errors and systematic errors. \citet{huan2013simulation, huan2014gradient, shahmohammadi2020using} and others have developed Bayesian MBDoE methods to consider different priors, error models, and incomplete data. Bayesian MBDoE formulates the objective function using information theoretic measures, such as the Kullback-Leibler (KL) divergence, to reflect the expected information gain.

Gradient-based methods are often applied to directly solve optimization problem \eqref{eq:simple-opt}. For example, the popular modeling platform \texttt{gPROMS} directly supports MBDoE using the \texttt{SRQPD} (sequential quadratic programming) solver \citep{galvanin2013general}. \citet{wang2022pyomo} proposed a two-stage program to formulate and solve Eq.~\eqref{eq:simple-opt} in their implementation of the open-source MBDoE software \texttt{Pyomo.DoE}. However, gradient-based solvers often get stuck at local optima. \citet{chachuat2004new, papamichail2004global, singer2006global, lin2006deterministic, bajaj2020global, kappatou2022global} studied
deterministic global optimization methods for parameter estimation, which is closely related to MBDoE. Surrogate models have been employed to approximate the objective function to accelerate MBDoE optimization \citep{paulson2019optimal, petsagkourakis2021safe} including Bayesian optimization \citep{shao2022active, shao2023preferenceguided}. Enumeration, which can be accelerated with the sensitivity analysis of nonlinear programming (NLP) problems \citep{thierry2019nonlinear}, enables visualization to develop intuition about the mathematical model and helps identify good initial points for computational optimization \citep{wang2022pyomo}.

Alternatively, \citet{kusumo2022risk} proposed the continuous-effort design by solving: 
\vspace{-5pt}
\begin{equation}
    \{p_1^*, ..., p_{N_s}^*\} = \underset{\mathbf{p}}{\text{argmax}} ~\Psi \left(\sum_{j=1}^{N_s} p_j \bigM(\bigtheta,\bigphi_j) \right), \quad \sum_{j=1}^{N_s} p_j = N_e
\label{eq:pydex-opt}
\end{equation}

Eq.~\eqref{eq:pydex-opt} discretizes the experimental design space into a set of $N_s$ candidate experiments \{$\bigphi_1, ..., \bigphi_{N_s}$\}, and decides the weights \{$p_1^*, ..., p_{N_s}^*$\} assigned to each candidate experiment by maximizing the FIM to select the best set of $N_e$ experiments. The FIM contributions $\bigM(\bigtheta,\bigphi_j)$ for experimental conditions $\bigphi_j$ are pre-computed, thus avoiding evaluating nonlinear dynamical models during optimization. This problem formulation is also concave for popular MBDoE metrics, including A-, D-, and E- optimality, regardless of the model structure \citep{vandenberghe1998determinant}. These capabilities are available in the software \texttt{Pydex} \citep{kusumo2022risk}. Specifically, \citet{kusumo2022risk} constructed the sensitivity equations by automatic differentiation using \texttt{CasADi} \citep{andersson2019casadi} interfaced through a modified version of \texttt{Pyomo.DAE} \citep{nicholson2018pyomo}, then computed the atomic matrices $\bigM(\bigtheta,\bigphi_j)$ using the forward sensitivity analysis capability of the solver \texttt{IDAS}. The convex formulation \eqref{eq:pydex-opt} with $p_i$ relaxed as continuous variables was solved using \texttt{MOSEK} interfaced through \texttt{CVXPY} which automatically applied a conic reformulation of the log determinant objective (D-optimality). However, there can be significant memory requirements for storing $\bigM(\bigtheta,\bigphi_j)$ for many candidate experiments $N_s$. Recently, \citet{hendrych2023solving} benchmarked different mixed integer solvers for Eq.~\eqref{eq:pydex-opt} and related experiment design problems using the package \texttt{Boscia.jl}. Our work pursues an alternate strategy to incorporate log determinant objectives (or constraints) into \texttt{Pyomo} using the \texttt{ExternalGreyBoxModel} module.

\subsection{Sensor and Measurement Selection Optimization} 

Sensor network design (SND) is closely related to selecting measurements in an experiment.
Most SND literature considered mixed-integer optimization to choose the best sensors for dynamic chemical processes, including power plants with CO$_2$ capture \citep{paul2015sensor, paul2016dynamic}, water treatment \citep{murray2010sensor}, water networks \citep{mann2012real, seth2016testing}, gas and flame detection \citep{benavides2016optimal, zhen2019mathematical},  and packed bed reactors \citep{serpas2013sensor}. 
Common objectives of SND include minimizing the sensor network costs \citep{bagajewicz2002new, chmielewski2002theory, kelly2008new}, maximizing estimation accuracy \citep{kadu2008optimal}, maximizing estimation reliability \citep{ali1993sensor, ali1995redundant}, and maximizing process robustness \citep{bhushan2008robust}. 
These objectives can also be incorporated as constraints within the SND framework, depending on the specific application and requirements \citep{bagajewicz2002new, paul2015sensor, paul2016dynamic, mobed2017state}.

Most applications of SND have focused on continuous measurements where a sensor can dynamically observe one measurement at many sampling times \citep{zhen2019optimal, legg2013optimal, benavides2016optimal}. In contrast, limited literature \citep{klise2017sensor, klise2020sensor, mann2012real} considered dynamic-cost measurements, such as manual sampling, in SND problems. \citet{klise2017sensor} developed an open-source Python package \texttt{Chama} for SND and applied it to a water network model. In this package, they classified sensors as stationary or mobile according to whether they are fixed at one location, and their capabilities to capture one or more signals. For stationary sensors, regardless of being a point or camera sensor, the amount and time points of samples need to be decided with a cost for each time point, which aligns with the concept of dynamic-cost measurements. \citet{mann2012real} considered manual grab samples providing a discrete indication of the presence of one component and proposed a mixed-inter linear programming (MILP) formulation to solve the problem using such discrete measurements at limited points in time and space. 

Several recent SND studies improved parameter estimation precision through FIM-based objectives based on local sensitivity information \citep{wouwer2000approach, qureshi1980optimum, basseville1987optimal, castro2013robustness, kretsovalis1987effect}. 
However, formulating these SND problems as a monolithic mixed-integer nonlinear programming (MINLP) problem is computationally difficult. 
\citet{wouwer2000approach} maximized an FIM-based criterion for a catalytic fixed bed reactor system, with partial differential-algebraic equations as constraints. They noted the high computational expense of this monolithic MINLP and instead trained
linear surrogate models. Similarly, \citet{isidori1985nonlinear, hermann1977nonlinear, lopez2004effect, yamada2022greedy} addressed nonlinear systems using data-driven linear reduced-order models. 
\citet{muske2003optimal} presented sensor location techniques to compromise between measurement costs and process information for parameter estimation and applied the technique to linear systems. 
\citet{serpas2013sensor} proposed a related problem that maximized the determinant of the observability matrix to optimize certain properties of the process state estimates. \citet{serpas2013sensor} presented a decomposition strategy to mitigate numerical issues with the determinant for ill-conditioned problems.

\subsection{Paper Contribution: Measurement Optimization}

This paper proposes a convex MI(N)LP formulation to select measurements for maximizing the information for validating science-based mathematical models. To ensure computational tractability, the D-optimality objective (i.e., determinant of the FIM) and its gradient are calculated in \texttt{SciPy} and incorporated into Pyomo using the new \texttt{ExternalGreyBoxModel} feature \citep{rodriguez2023scalable}. A modified version of \texttt{MindtPy} is developed to support grey-box constraints.  The scalability of our approach is demonstrated in two case studies. Case study 1 considers a nonlinear reaction kinetics model with three highly correlated time-varying measurements, which shows how MO balances between budgets and information contents.  Case study 2 considers a large-scale rotary packed bed for CO$_2$ capture, with 14 time-varying measurements that are measurable at hundreds of discretized time steps, formulating an optimization problem of millions of variables.

\section{Methodology}

An experiment shown in Eq.~\eqref{eq:v1} consists of specifying experimental conditions $\bigphi \in \boldsymbol{\Phi} \subset \R^{C}$ and measuring a subset of the state variables, i.e. dynamical measurements, $\bigy_k\in \mathcal{Y}  \subset \R^{K} =\{\bigy_1,...,\bigy_{K}\}$ at a set of specific time points $\bigt_k=\{t_1,...,t_{T_{k}}\}$, to estimate uncertain parameters $\bigtheta \in \boldsymbol{\Theta} \subset \R^{P}$. Thus, $y_{k,t}$ denotes the measurement $y_k$ at a specific time point $t$. The error covariance matrix $\bigSigma_{\epsilon}$ corresponds to the observation error $\boldsymbol{\epsilon}$ in Eq.~\eqref{eq:v1}, which corrupts all measured data. 

Dynamic-cost measurements (DCMs) are denoted by $\mathcal{D}=\{\bigy_d: \forall d\in \{1,...,D\}\}$. DCMs, referred to as manual grab samples in SND literature \citep{mann2012real}, have a fixed installation cost $c_k$ and a significant per-measurement cost $c_{k,t}$, for example, composition analysis that requires manual measurements at specific times by a technician. DCMs are often constrained by a minimum time between each measurement or a maximum of total measurements in a time period. The binary variable $x_{k,t}\in\{0,1\}, \forall d \in \mathcal{D}$ indicates if the data $y_{k,t}$ is measured at time point $t$, with a cost $c_{k,t}$. 

Static-cost measurements (SCMs) are a special case of DCMs, which have a one-time cost to install the sensor and a negligible edge cost per measurement. SCMs are denoted by $\mathcal{S}=\{\bigy_s, \forall s\in\{1,...,S\}\}$. SCMs are ``all or nothing,'' i.e., if the sensor is installed, all of the corresponding time points are measured. Thermocouples and inline gas chromatography (GC) machines are SCMs since, once selected and installed at a fixed cost, their values at various timestamps through the experiment can be measured with no or minimal marginal costs. The binary variable $x_k\in\{0,1\}, \forall k \in \mathcal{S}$ is not indexed by time for SCMs.

\subsection{Dynamic Sensitivity Analysis}

For dynamic models, the measured variable $\bigy_k$ is sampled at $\bigt_k=\{t_1,...,t_{T_{k}}\}$. All of the measurements can be collected in a single vector $\bigy$: 
\vspace{-5pt}
\begin{equation}
    \bigy = [y_{1,t_1}, y_{1,t_2}, ..., y_{K, t_{T_K}}]
\label{eq:ystack}
\end{equation}

\noindent where $y_{k,t}$ represents measurement $k$ at time $t$. 

$\bigQ$, the dynamic sensitivity matrix for all the measurements $\bigy$, is defined as:

\vspace{-5pt}
\begin{equation}
    \bigQ = \nabla_{\bigtheta}~\bigy
\label{eq:Q}
\end{equation}

\noindent $\bigQ$ includes the sensitivity of all measurements. Scalar  $q_{(d,t),i} \in \R^{1\times1}$ is the dynamic sensitivity $\frac{\partial y_{d,t}}{\partial \theta_i}$ for a DCM $y_{d,t}, d\in\mathcal{D}$ at time $t\in\mathbf{t}_d$ with respect to parameter $\theta_i$. Vector $\bigq_{s,i} \in \R^{T_s \times 1}$ is the dynamic sensitivity $\frac{\partial y_{s}}{\partial \theta_i}$ for SCM $y_s, s\in\mathcal{S}$ at all time points $\bigt_s$ with respect to parameter $\theta_i$.

\subsection{Information Content of Measurements}

The FIM $\bigM \in \R^{P\times P}$ quantifies the information about unknown parameters $\bigtheta=\{\theta_1, ..., \theta_P\}$ encoded in measurements $\bigy$. The parameter covariance matrix $\mathbf{V}$ quantifies the uncertainty in the estimated parameter values $\hat{\bigtheta}$:
\vspace{-5pt}
\begin{equation}
    \mathbf{V} \approx (\bigQ^\intercal\bigSigma_{\epsilon}^{-1}\bigQ)^{-1}
\label{eq:v}
\end{equation}

\noindent where $\bigSigma_{\epsilon}^{-1}$ is the inverse error covariance matrix and $\bigQ$ is the full dynamic sensitivity matrix. 
The error covariance matrix $\bigSigma_{\epsilon}$ captures all the pairwise correlations between measurements across different time steps. The modeler must choose the appropriate structure for $\bigSigma_{\epsilon}$. If the modeler believes the measurement errors are truly independent, then $\bigSigma_{\epsilon}$ is a diagonal matrix. Otherwise, the modeler can fit an error structure model, e.g., autoregressive or Gaussian Process \citep{lainez2015stochastic}.

The FIM can be estimated as the inverse of $\mathbf{V}$ using $\bigSigma_{\epsilon}$ and $\bigQ$:
\vspace{-5pt}
\begin{equation}
    \bigM \approx \mathbf{V}^{-1} \approx \bigQ^\intercal\bigSigma_{\epsilon}^{-1}\bigQ
\label{eq:mv}
\end{equation}

\noindent Eqs.~\eqref{eq:v} and \eqref{eq:mv} asymptotically hold for nonlinear models \citep{bard1974nonlinear}. We now consider how to compute the contributions of SCMs and DCMs to $\bigM$.

Let $\bigM_{s, s'}$ be the FIM of two SCMs $s$ and $s'$. Its element 
at the $i^{\mathrm{th}}$ row, $j^{\mathrm{th}}$ column is computed by: 
\vspace{-5pt}
\begin{equation}
    m_{s, s', i, j} \approx \bigq_{s,i}^\intercal\cdot\tilde{\bigSigma}_{s,s'}\cdot\bigq_{s',j}, \forall i, j \in \{1,...,P\}
\end{equation}

\noindent where $\bigq_{s,i} \in \R^{T_s\times1}$, $\bigq_{s',i} \in \R^{T_{s'}\times 1}$, $\tilde{\bigSigma}_{s,s'} \in \R^{T_s\times T_{s'}}$ are the elements of the inverse of the error covariance matrix $\bigSigma_{\epsilon}^{-1}$ corresponding to the two SCMs.

$\bigM_{(d,t), (d',t')}$ is the FIM of the DCM $d$ at time point $t$, and the DCM $d'$ at time point $t'$. Its element at the $i^{\mathrm{th}}$ row, $j^{\mathrm{th}}$ column is computed by:
\vspace{-5pt}
\begin{equation}
    m_{(d,t), (d',t'), i, j} \approx q_{(d,t),i}\cdot\tilde{\sigma}_{(d,t),(d',t')}\cdot q_{(d',t'),j}, \forall i,j \in \{1,...,P\},
\end{equation}

\noindent where $q_{(d,t),i} \in \R^{1 \times 1}$, $q_{(d',t'),j} \in \R^{1 \times 1}$, and $\tilde{\sigma}_{(d,t),(d',t)} \in \R^{1\times1}$ are the elements of the inverse of the error covariance matrix $\bigSigma_{\epsilon}^{-1}$ corresponding to the two DCMs. 

$\bigM_{s, (d,t)}$ is the FIM of the DCM $d$ at the time point $t$  and the SCM $s$. Its element at the $i^{\mathrm{th}}$ row, $j^{\mathrm{th}}$ column is computed by:

\vspace{-5pt}
\begin{equation}
    m_{s, (d,t), i, j} \approx
  \bigq_{s,i}^\intercal\cdot\tilde{\bigSigma}_{s,(d,t)}\cdot q_{(d,t),j}, \forall i,j \in \{1,...,P\}
\end{equation}

\noindent where $\bigq_{s,i} \in \R^{T_s\times 1}$,  $q_{(d,t),j} \in \R^{1 \times 1}$, and $\tilde{\bigSigma}_{s,(d,t)} \in \R^{T_s\times 1}$ are the elements of the inverse of the error covariance matrix $\bigSigma_{\epsilon}^{-1}$ corresponding to the DCM and the SCM. 

\subsection{Measurement Optimization Framework}

We propose the following MO problem to select the best set of measurements under a specific budget and manual sampling constraints:
\vspace{-8pt}
\begin{subequations}
\begin{align}
\underset{\mathbf{x}}{\max} \ \ \  &\Psi(\bigM + \bigM_0) \\
\text{s.t.} \ \ \ m_{i,j} = &\sum_{d\in\mathcal{D}}\sum_{d'\in\mathcal{D}}\sum_{t\in\mathbf{t}_d} \sum_{t'\in\mathbf{t}_{d'}} m_{(d,t), (d',t), i, j}\cdot x_{(d,t), (d',t')} \label{eq:dynamic-M}\\ 
&+\sum_{s\in\mathcal{S}}\sum_{d\in\mathcal{D}} \sum_{t\in\mathbf{t}_{d}} m_{s, (d,t), i, j}\cdot x_{s, (d,t)}\label{eq:static-dynamic-M}\\ &+\sum_{s\in\mathcal{S}}\sum_{s'\in\mathcal{S}} m_{s, s', i, j}\cdot x_{s, s'}, \forall i,j\in\{1,...,P\} \label{eq:static-M} \\
& \sum_{k\in \mathcal{Y}} c_k \cdot x_k \leq B \label{eq:con_cost}\\
& \sum_{d\in\mathcal{D}}\sum_{t\in\bigt_d} x_{d,t} \leq L_{total} \label{eq:ltotal}\\
& \sum_{t\in\bigt_d} x_{d,t} \leq L_d, \ \ \forall d \in \mathcal{D} \label{eq:ld} \\ 
& \sum_{t'}^{t'-t<T_{min}}x_{d,t} \leq 1, \ \ \forall d \in \mathcal{D}, \forall t \in \bigt_d \label{eq:interval} \\
& x_{k,k'} \leq x_{k,k}, \ \ \forall k < k', k, k' \in \mathcal{Y} \label{eq:con1} \\ 
& x_{k,k'} \leq x_{k',k'}, \ \ \forall k < k', k, k' \in \mathcal{Y} \\ 
 &  x_{k,k}+x_{k',k'}-1 \leq x_{k,k'}, \ \ \forall k<k', k, k' \in \mathcal{Y} \label{eq:con3}
\end{align}
\label{eq:overall}
\end{subequations}

\vspace{-15pt}
\begin{table}[htb!]
\centering
\caption{Measurement optimization problem size.} 
\begin{tabular}{lc}
\toprule
& Problem size\\
\midrule
Variables & $(N_r^2+N_r)/2+(P^2+P)/2 + D + 2 $) \\ 
Equality Constraints & $(P^2+P)/2  +2$ \\
Inequality Constraints & $3(N_r^2-N_r)/2 + 2 + 2D + 2\Sigma_d^D T_d$ \\
\bottomrule
\end{tabular}
\label{tab:size}
\end{table}

The binary measurement selection variables $\mathbf{x}=\{x_{1}, ..., x_{{S}}, x_{{1,t_1}},..., x_{{D, t_{T_D}}}\}$ are manipulated to maximize the information content of experiments. $\bigM_0$ is the prior FIM.
Eq.~\eqref{eq:dynamic-M} to Eq.~\eqref{eq:static-M} calculate the total FIM by summing the contributions from the selected measurements. Eq.~\eqref{eq:con_cost} constrains the total cost of measurements to be less than the total budget $B$ considering the cost $c_k$ of measurement $k\in \mathcal{Y}$. Eq.~\eqref{eq:ltotal} constrains the maximum total number of measurements for all dynamic-cost measurements to be less than or equal to $L_{total}$. Eq.~\eqref{eq:ld} limits the maximum total number of measurements for one dynamic-cost measurement $d$ to be less than or equal to $L_d$. Eq.~\eqref{eq:interval} requires a minimal time interval $T_{min}$ between two dynamic-cost measurements. Parameters $L_{total}$, $L_d$, and $T_{min}$ are specified by the modeler. Table \ref{tab:size} reports the size of MO problem \eqref{eq:overall} as a function of the total number of measurements $N_r=S+\sum_{k\in\mathcal{D}}T_k$, $P, D$, and $T_D$.

Recall that the binary decision $x_k\in\{0,1\}$ indicates if measurement $k$ is selected. Likewise, the binary decision $x_{k,k'} \in \{0,1\}$ indicates if measurements $k$ and $k'$ are both selected. Eq.~\eqref{eq:con1} to Eq.~\eqref{eq:con3} are the McCormick relaxations for the bilinear equation $x_{k,k}\cdot x_{k',k'} = x_{k,k'}$. They are used to ensure the covariances between two measurements are included if and only if both measurements are selected. 
All constraints in optimization problem \eqref{eq:overall} are linear in $\mathbf{x}$. Relaxing $x_k\in\{0,1\}$ to $x_k\in [0,1]$ results in a continuous linear programming (LP) problem.  

\subsection{Computational Implementation}

Optimization problem \eqref{eq:overall} is implemented in \texttt{Pyomo} version 6.7 \citep{bynum2021pyomo}. Two design criteria, A-optimality (trace) and D-optimality (determinant), are separate objective functions. A-optimality computes the trace by summing up the diagonal elements of $\bigM$:
\vspace{-5pt}
\begin{equation}
   \trace(\bigM) = \sum_{i=1}^{P}m_{i,i} 
\end{equation}

\noindent The MILP problem and its relaxed LP problem are solved using \texttt{Gurobi} version 10.0.3. 

D-optimality, $\det(\bigM)$, is a concave function, constituting a convex maximization problem. The objective function is often formed as: 
\vspace{-5pt}
\begin{equation}
    \psi(\bigM) = \log(\det(\bigM))
\label{eq:obj}
\end{equation}

\noindent However, computing determinants algebraically is challenging. \citet{wang2022pyomo} used the Cholesky factorization $\bigM=\mathbf{L}\mathbf{L}^\intercal$, where $\mathbf{L}$ is a lower triangular matrix, converting the determinant calculation into a summation of bilinear terms. The objective function \eqref{eq:obj} is evaluated with the diagonal elements of $\mathbf{L}$ by: 
\vspace{-5pt}
\begin{equation}
    \log(\det(\bigM))=2\sum_{i=1}^P \log (l_{ii})
\end{equation}

\noindent Unfortunately, adding the Cholesky factorization as an equality constraint introduces bilinear expressions, making the optimization problem non-convex. Alternately, \citet{kusumo2022risk} and \citet{hendrych2023solving} consider conic reformulations of the log determinant objective and use convex optimization solvers. This approach works well for the continuous effort formulation \eqref{eq:pydex-opt} where the relaxed problem is convex. However, the direct MBDoE problem \eqref{eq:simple-opt}, where $\phi$ is directly optimized (instead of the decision space being discretized \emph{a priori}), is often non-convex due to nonlinear model constraints.

This motivates our alternate strategy to calculate log determinants in an algebraic modeling language, e.g., \texttt{Pyomo}, and use standard non-convex nonlinear programming solvers, e.g., \texttt{Ipopt} \citep{wachter2006implementation}. Specifically, we use the \texttt{ExternalGreyBoxModel} feature in \texttt{Pyomo} to compute the objective function \eqref{eq:obj} and its gradient \citep{de2006aspects} with \texttt{SciPy}: 

\vspace{-5pt}
\begin{equation}
    \nabla_{\bigM} \psi(\bigM) = \bigM^{-1}
\label{eq:greybox1}
\end{equation}

Because the FIM $\bigM$ is symmetric, real-valued, and positive-semi-definite, Eq.~\eqref{eq:greybox1} simplifies as follows \citep{harville1998matrix}: 

\vspace{-5pt}
\begin{equation}
    \nabla_{\bigM} \psi(\bigM) = 2\bigM^{-1} - \text{diag}(\bigM^{-1})
    \label{eq:grey}
\end{equation}

For computation, we exploit the symmetry of $\bigM$ by defining the vector $\mathbf{m}=[ m_{11}, m_{12}, ..., m_{1P}, m_{22}, ..., m_{PP} ]^\intercal$ as the lower triangular elements of $\bigM$. We evaluate Eq.~\eqref{eq:grey} using the Moore-Penrose pseudo-inverse in \texttt{SciPy} but return $\nabla_{\textbf{m}} \psi(\bigM) = \left[ \frac{\partial \psi(\bigM)}{\partial m_{11}}, ..., \frac{\partial \psi(\bigM)}{\partial m_{PP}} \right]^\intercal$ to interface with \texttt{ExternalGreyBoxModel} in \texttt{Pyomo}. The pseudo-inverse provides some robustness if $\bigM$ is ill-conditioned. If $\bigM$ is rank deficient, we recommend using a sufficiently large prior FIM $\bigM_0$ as a regularization term.

The optimization problem \eqref{eq:overall} becomes an MINLP problem with the nonlinear objective function \eqref{eq:obj}. A modified version of \texttt{MindtPy} version 1.0.0 \citep{bernal2018mixed} with a customized warm-start initialization was used to solve the D-optimality MINLPs through the outer approximation method using \texttt{Gurobi} version 10.0.3 for the MILP subproblems. \texttt{CyIpopt} version 1.3.0 \citep{cyipopt} , a C-extensions for Python interface to \texttt{Ipopt} version 3.14.12, and \texttt{MA57} \citep{hsl2018collection} were used with the Pyomo \texttt{ExternalGreyBoxModel} module to solve the NLP subproblems. 

For both case studies, we solved the optimization problems with the lowest budget and gradually increased to the highest budget using the prior solution as an initial point. We used the A-optimality LP solutions to initialize the MILP problems. We then used these solutions to initialize the D-optimality NLP and MINLP problems. Moreover, we developed a customized strategy to warm-start the fixed NLP subproblems for each iteration. Following the solution of the MILP subproblem, our customized warm-start strategy explicitly calculated the continuous variable values in the model based on the MILP integer solution using \texttt{SciPy}. These values were then used to initial the continuous variables in the fixed NLP subproblem, thus provided a feasible starting point. Unless otherwise specified, all reported computational times refer to the solver only and do not include model building in \texttt{Pyomo}. All the experiments were conducted on dual 12-core Intel Xeon CPU E5-2680 with 256 GB of RAM.

For validation, optimization problem \eqref{eq:overall} was also implemented in \texttt{CVXPY} 1.2.1 \citep{diamond2016cvxpy, agrawal2018rewriting}, a Python-based modeling environment for convex optimization problems. We used \texttt{Mosek} 9.3.20 \citep{mosek} to solve the A-optimality MILP and D-optimality NLP problems and confirmed similar solutions to our proposed approach in \texttt{Pyomo}. At the time of writing, the \texttt{CVXPY} interface to \texttt{Mosek} did not support mixed integer optimization. These comparison results were generated with a Macbook Pro (15-inch, 2018) with a 2.6 GHz 6-core Intel Core i7 processor and 16 GB of RAM.

\section{Case Study: Reaction Kinetics}

We start by illustrating the MO framework using a reaction kinetics example with highly correlated measurements. 

\subsection{Problem Statement}

Consider two first-order liquid phase reactions in a batch reactor:

\begin{equation*}
    \text{A} \xrightarrow{k_1} \text{B} \xrightarrow{k_2} \text{C}
\end{equation*}

\begin{table}[htb!]
\centering
\caption{Mathematical model for the reaction kinetics example. The Arrhenius equations model the temperature dependence of the reaction rate coefficients $k_1$ and $k_2$. We assume a first-order reaction mechanism and only species A is fed to the reactor.}  
\begin{tabular}{lc}
\toprule
&Equations\\
\midrule
Arrhenius equation & $k_1 = A_1 e^{-\frac{E_1}{RT}} $\\
&$k_2 = A_2 e^{-\frac{E_2}{RT}}$\\
\\
Reaction rates & $\frac{d{C_A}}{dt} = -k_1{C_A}$ \\
& $\frac{d{C_B}}{dt} = k_1{C_A} - k_2{C_B}$ \\
\\
Mole balance & $C_A + C_B + C_C = C_{A0}$ \\
\\
Initial conditions & $C_A(t_0) = C_{A0}$ \\ 
& $C_B(t_0) = 0$ \\
& $C_C(t_0) = 0$ \\ 
\bottomrule
\end{tabular}
\label{tab:reactor}
\end{table}

\begin{figure}[!ht]
\centering
\includegraphics[width=0.6\textwidth]{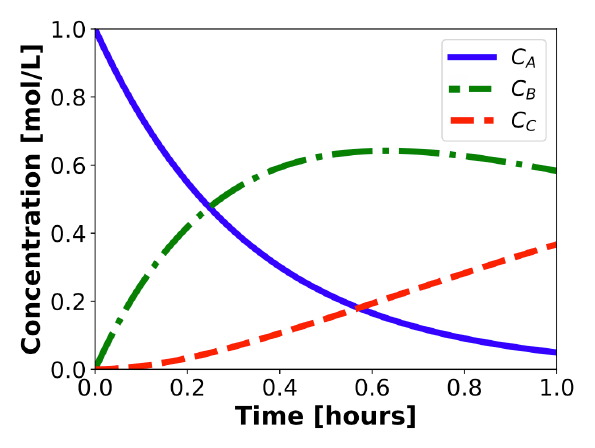}
\caption{Time-varying profiles for the concentration of component A (blue line), B (green dash line), C (red dash line).} 
\label{fig:kinetics-setup}
\end{figure}

\noindent Table \ref{tab:reactor} gives the mathematical model for the system \citep{wang2022pyomo}. Measurements $\mathbf{y} = [\mathbf{C_A}, \mathbf{C_B}, \mathbf{C_C}]^\intercal$ are the time-varying concentrations of the species A, B, and C. Unknown parameters are the activation energies $E_1$ and $E_2$, and pre-exponential factors $A_1$ and $A_2$ in the Arrhenius equations computing the reaction rate constants $k_1$ and $k_2$. Time-varying profiles for the concentration of all species are shown in Fig.~\ref{fig:kinetics-setup}. MO determines the measurement strategy to maximize the precision of the estimates of model parameters $\bigtheta$=[$A_1, E_1, A_2, E_2$]$^\intercal$.

\begin{table}[htb!]
\centering
\caption{Assumed costs for SCMs and DCMs for the kinetics case study.}  
\begin{tabular}{|l |c |c |}
\hline 
Name & Installation cost (\$) & Measurement cost (\$/sample) \\ \hline
$c_A^{SCM}$ & 2000 & 0 \\
$c_B^{SCM}$ & 2000 & 0 \\
$c_C^{SCM}$ & 2000 & 0 \\ \hline
$c_A^{DCM}$ & 200 & 400 \\
$c_B^{DCM}$ & 200 & 400 \\
$c_C^{DCM}$ & 200 & 400 \\ \hline

\end{tabular}
\label{tab:kinetics_price}
\end{table}

We assume $C_A$, $C_B$, and $C_C$ are measured via high-performance liquid chromatography (HPLC). HPLC machines with automated sample extraction equipment can be integrated directly into the experimental system to measure concentrations continuously and thus are SCMs; HPLC machines with manual sampling extraction are modeled as DCMs. Table \ref{tab:kinetics_price} lists the measurement costs based on \citet{liptak2003instrument}. Each time-varying measurement $C_i$ is discretized to nine time points evenly from 0 to 60 minutes.  Considering human resources constraints, a maximum of five time points can be chosen for each DCM, and ten are allowed for all DCMs. Additionally, the selected time points must maintain a minimum separation of at least ten minutes between each other, which means there should be no two consecutive time points for DCMs under this time discretization. For completeness, the Supporting Information (SI) also analyzes a minor operating cost of the SCMs (section \ref{sec:operate}) and HPLCs with multi-component samples (section \ref{sec:multi}). The main findings of these supplemental analyses are the same as the case considered in the main text.

The observation errors are assumed to be independent across different time points, but with correlations between measurements:
\vspace{-5pt}
\begin{equation}
\boldsymbol{\Sigma_{\epsilon}}=\left[\begin{array}{ccc}
\Var(C_A) & \Cov(C_A, C_B) & \Cov(C_A, C_C) \\
 \Cov(C_A, C_B) & \Var(C_B) & \Cov(C_B, C_C) \\
\Cov(C_A, C_C) & \Cov(C_B, C_C) & \Var(C_C) \\
\end{array}\right]=\left[\begin{array}{ccc}
1 & 0.1 & 0.1 \\
0.1 & 4 & 0.5 \\
0.1 & 0.5 & 8 \\
\end{array}\right]
\label{eq:Cov-kinetics}
\end{equation}

\noindent where each (co)variance term units of mol$^2 \cdot \text{L}^{-2}$. The covariance between one DCM and one SCM is half of its original correlation value. SI Table \ref{tab:kinetics_cov} reports the full covariance structure.

Four optimization problems are solved: maximizing A-optimality of the MILP problem; its relaxed LP problem; maximizing D-optimality of the MINLP problem; and its relaxed NLP problem.
The A-optimality LP problem has 393 continuous variables, 12 equality constraints, and 1,096 inequality constraints. 
The A-optimality MILP problem has 12 continuous variables, 381 binary variables, 12 equality constraints, and 1,096 inequality constraints. 
The D-optimality NLP problem has 404 variables, 23 equality constraints, and 1,096 inequality constraints. 
The D-optimality MINLP problem has 23 continuous variables, 381 binary variables, 23 equality constraints, and 1,096 inequality constraints. 

\setlength\tabcolsep{3pt} 
\begin{longtable}{|l | c | c | c | c | c | c| c| c |c |c |c |}
\caption{Optimal solutions for the MILP A-optimality optimization problem where each column corresponds to a measurement budget of \$1 k to \$5 k. '1' and '0' indicate if each SCM is selected or not selected.  The selected sample times [min] of DCMs are reported.}\\
\hline
\multicolumn{1}{|c|}{\textbf{Budget [\$ k]}} & 
\multicolumn{1}{c|}{\textbf{1.0}} & 
\multicolumn{1}{c|}{\textbf{1.4}} &
\multicolumn{1}{c|}{\textbf{1.8}} &
\multicolumn{1}{c|}{\textbf{2.2}} &
\multicolumn{1}{c|}{\textbf{2.6}} &
\multicolumn{1}{c|}{\textbf{3.0}} &
\multicolumn{1}{c|}{\textbf{3.4}} &
\multicolumn{1}{c|}{\textbf{3.8}} &
\multicolumn{1}{c|}{\textbf{4.2}} &
\multicolumn{1}{c|}{\textbf{4.6}} &
\multicolumn{1}{c|}{\textbf{5.0}}\\ \hline 
\endfirsthead

\multicolumn{12}{c}%
{{\bfseries \tablename\ \thetable{} -- continued from previous page}} \\
\hline 
\multicolumn{1}{|c|}{\textbf{Budget [\$ k]}} & 
\multicolumn{1}{c|}{\textbf{1.0}} & 
\multicolumn{1}{c|}{\textbf{1.4}} &
\multicolumn{1}{c|}{\textbf{1.8}} &
\multicolumn{1}{c|}{\textbf{2.2}} &
\multicolumn{1}{c|}{\textbf{2.6}} &
\multicolumn{1}{c|}{\textbf{3.0}} &
\multicolumn{1}{c|}{\textbf{3.4}} &
\multicolumn{1}{c|}{\textbf{3.8}} &
\multicolumn{1}{c|}{\textbf{4.2}} &
\multicolumn{1}{c|}{\textbf{4.6}} &
\multicolumn{1}{c|}{\textbf{5.0}}\\ \hline 
\endhead

\hline \multicolumn{12}{|r|}{{Continued on next page}} \\ \hline
\endfoot

\hline \hline
\endlastfoot
 
$C_{A}^{SCM}$ & 0  & 0  & 0  & 0 
& 0& 0& 0 &0 
& 0& 0& 0 \\ 
$C_{B}^{SCM}$ & 0  & 0  & 0  & 1 
& 1& 1& 1 &1 
& 1& 1& 1 \\
$C_{C}^{SCM}$ & 0  & 0  & 0  & 0 
& 0& 0& 0 &0 
& 1& 1& 1  \\
\hline 
$C_{A}^{DCM}$ &   &   &  &  
& 7.5 & 7.5 &  &  
& & 7.5& 7.5  \\ 
&   &   &  &  
&  & 22.5 &  &  
& & & 22.5  \\

$C_{B}^{DCM}$ & 45  & 30  & 15 &  
&  &  &  &  
& & & \\ 
&  60 & 45  & 30 &  
&  &  &  &  
& & & \\ 
&   & 60  & 45  &  
&  &  &  &  
& & & \\ 
&   &   &  60 &  
&  &  &  &  
& & & \\ 
$C_{C}^{DCM}$ &   &  & &  
&  & &30  &15    
& & & \\ 
&   &   &  &  
&  & &45  &30    
& & & \\ 
&   &   &  &  
&  & &60   &45    
& & & \\ 
&   &   &   &  
&  &  & & 60   
& & & 
\label{tab:kinetics_solution}
\end{longtable}

\setlength\tabcolsep{3pt} 
\begin{longtable}{|l | c | c | c | c | c | c| c| c |c |c |c |}
\caption{Optimal solutions for the D-optimality MINLPs where each column corresponds to a measurement budget of \$1 k to \$5 k. '1' and '0' indicate if each SCM is selected or not selected.  The selected sample times [min] of DCMs are reported.}\\
\hline
\multicolumn{1}{|c|}{\textbf{Budget [\$ k]}} & 
\multicolumn{1}{c|}{\textbf{1.0}} & 
\multicolumn{1}{c|}{\textbf{1.4}} &
\multicolumn{1}{c|}{\textbf{1.8}} &
\multicolumn{1}{c|}{\textbf{2.2}} &
\multicolumn{1}{c|}{\textbf{2.6}} &
\multicolumn{1}{c|}{\textbf{3.0}} &
\multicolumn{1}{c|}{\textbf{3.4}} &
\multicolumn{1}{c|}{\textbf{3.8}} &
\multicolumn{1}{c|}{\textbf{4.2}} &
\multicolumn{1}{c|}{\textbf{4.6}} &
\multicolumn{1}{c|}{\textbf{5.0}}\\ \hline 
\endfirsthead

\multicolumn{12}{c}%
{{\bfseries \tablename\ \thetable{} -- continued from previous page}} \\
\hline 
\multicolumn{1}{|c|}{\textbf{Budget [\$ k]}} & 
\multicolumn{1}{c|}{\textbf{1.0}} & 
\multicolumn{1}{c|}{\textbf{1.4}} &
\multicolumn{1}{c|}{\textbf{1.8}} &
\multicolumn{1}{c|}{\textbf{2.2}} &
\multicolumn{1}{c|}{\textbf{2.6}} &
\multicolumn{1}{c|}{\textbf{3.0}} &
\multicolumn{1}{c|}{\textbf{3.4}} &
\multicolumn{1}{c|}{\textbf{3.8}} &
\multicolumn{1}{c|}{\textbf{4.2}} &
\multicolumn{1}{c|}{\textbf{4.6}} &
\multicolumn{1}{c|}{\textbf{5.0}}\\ \hline 
\endhead

\hline \multicolumn{12}{|r|}{{Continued on next page}} \\ \hline
\endfoot

\hline \hline
\endlastfoot
 
$C_{A}^{SCM}$ & 0  & 0  & 0  & 0 
& 0& 0& 0 &0 
& 1& 1& 1 \\ 
$C_{B}^{SCM}$ & 0  & 0  & 0  & 0 
& 1& 1& 1 &1 
& 1& 1& 1 \\
$C_{C}^{SCM}$ & 0  & 0  & 0  & 0 
& 0& 0& 0 &0 
& 0& 0& 0  \\
\hline 
$C_{A}^{DCM}$ &   &   & 7.5  & 7.5 
& 7.5 & 7.5 & 7.5  &  7.5
& & &   \\ 
&   &   & 22.5  & 37.5  
&  & 22.5 & 22.5 &  22.5
& & &   \\ 
&   &   &   & 
& &  & 37.5 & 37.5 
& & &   \\ 
&   &   &   & 
& &  &  &  52.5
& & &   \\ 

$C_{B}^{DCM}$ & 7.5  & 7.5  & 60 & 22.5
& &  &  &  
& & & \\ 
&  60 & 22.5  &  &  60
&  &  &  &  
& & & \\ 
&   & 60  &   &  
&  &  &  &  
& & & \\ 

$C_{C}^{DCM}$ &   &  & &  
&  & &  &  
& &7.5 &7.5 \\ 
&   &   &  &  
&  & & & 
& & &60 \label{tab:kinetics_solution_D}
\end{longtable}

\begin{figure}[!ht]
\centering
\includegraphics[width=0.95\textwidth]{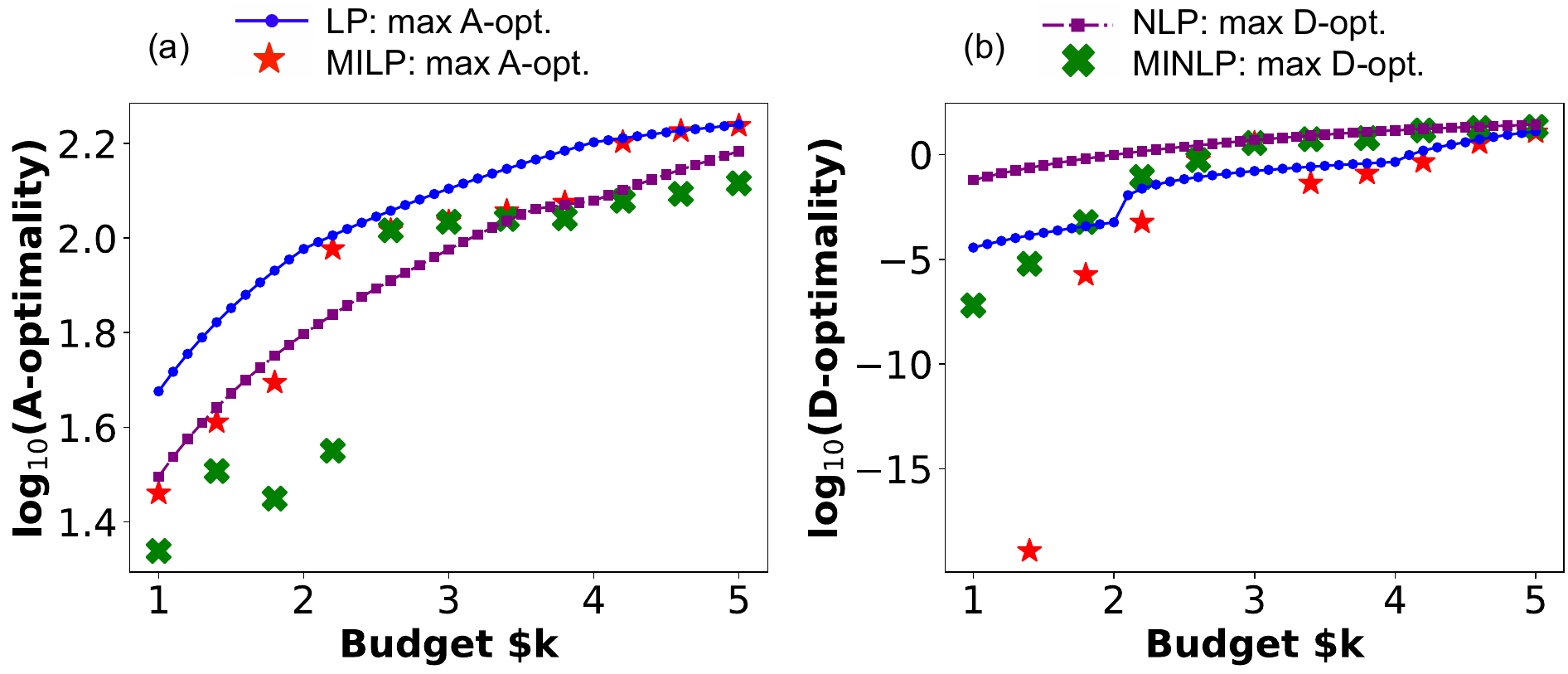}
\caption{Pareto-optimal trade-off between measurement budgets versus (a) A-optimality (trace of FIM) and (b) D-optimality (determinant of FIM) for the kinetics case study considering four optimization strategies: maximizing A-optimality of the MILP problem (red stars), its relaxed LP problem (blue line); maximizing D-optimality of the MINLP problem (green crosses), and its relaxed NLP problem (purple line). In (a), the blue line is an upper bound for the red stars, while in (b) the purple line is an upper bound for the green crosses.} 
\label{fig:kinetics}
\end{figure}

\begin{figure}[!ht]
\centering
\includegraphics[width=\textwidth]{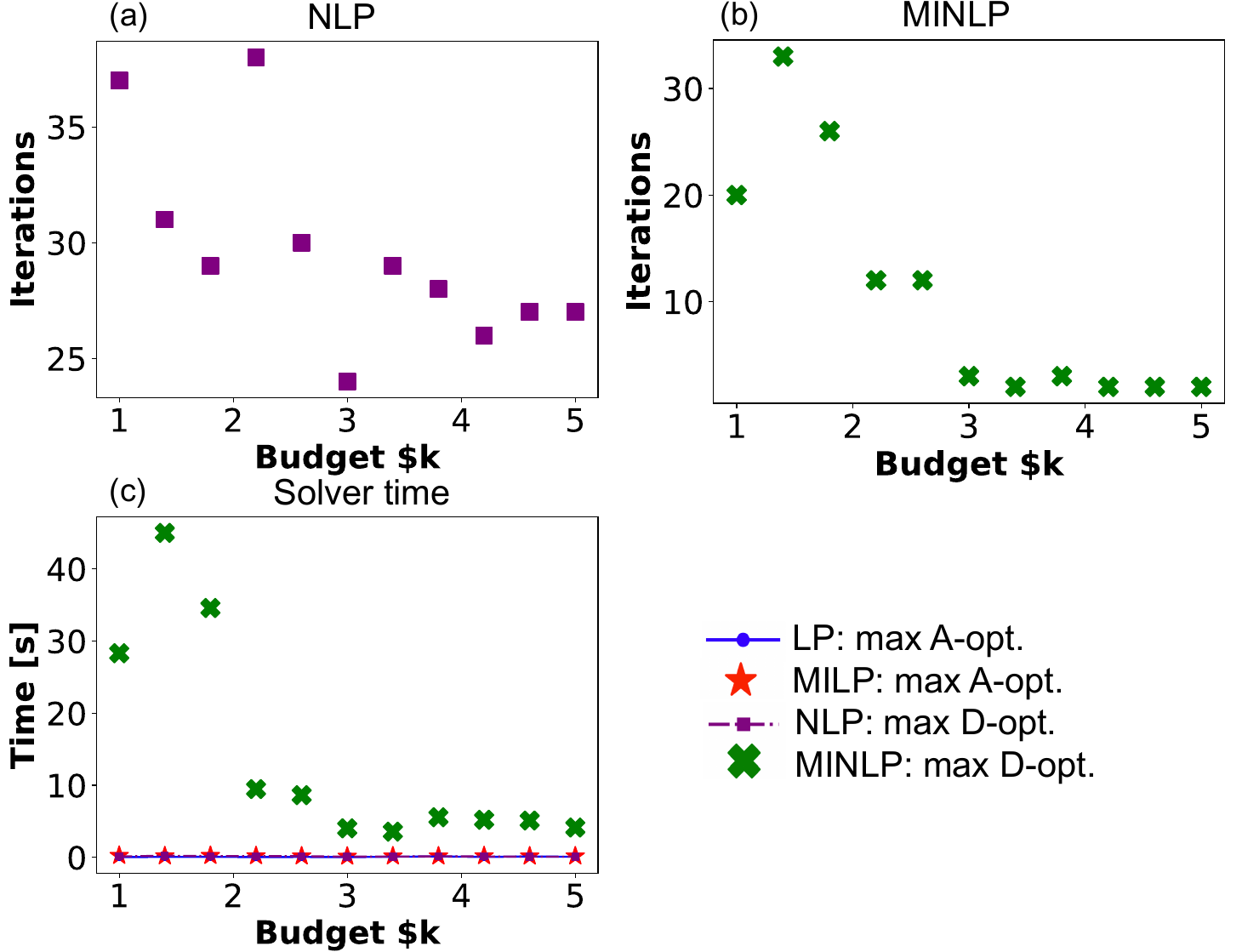}
\caption{Computational results for the kinetic case study at different budgets. (a) the number of \texttt{CyIpopt} iterations of D-optimality NLPs, (b) the number of \texttt{MindtPy} iterations of the D-optimality MINLPs, and (c) the computational time for all four optimization strategies. The optimization problems are solved from the lowest budget to the highest budget.} 
\label{fig:kinetics-calc}
\end{figure}

\subsection{Pareto-Optimal Trade-offs Between Measurement Budgets and Information}

Fig.~\ref{fig:kinetics} shows the Pareto-optimal trade-offs between measurement costs and information content, i.e., A-optimality and D-optimality metrics, and allows practitioners to select the best measurement budget when designing an experimental apparatus. A solution is Pareto-optimal if improving one objective, e.g., the information content, requires sacrificing one or more other objectives, e.g., the budget \citep{dowling2016framework}. Fig.~\ref{fig:kinetics} was generated by solving four optimization problems -- the MILP problem maximizing A-optimality, its relaxed LP problem, the MINLP problem maximizing D-optimality, and its relaxed NLP problem -- for 11 budgets from \$1 k to \$5 k. This is known as the $\epsilon$-constrained method for multi-objective optimization. Fig.~\ref{fig:kinetics}(a) shows the A-optimality values for the solutions for all four optimization problems, while Fig.~\ref{fig:kinetics}(b) shows the D-optimality values for these solutions. In Fig.~\ref{fig:kinetics}(a), the blue line is an upper bound for the red stars, while in Fig.~\ref{fig:kinetics}(b) the purple line is an upper bound for the green crosses.
Both A- and D-optimal solutions show substantial improvements as the budget increases from \$1 k to \$5 k. 
As anticipated, A-optimal solutions (i.e., maximize A-optimality) have higher A-optimality values than D-optimal solutions (i.e., maximize D-optimality). Likewise, D-optimal solutions show higher D-optimality value than A-optimal solutions. Fig.~\ref{fig:kinetics}(a) shows a sharp increase in A-optimality (red stars) at \$2.2 k with minimal increase until reaching a budget of \$4.2 k. Thus, a practitioner may conclude that A-optimality recommends \$2.2 k or \$4.2 k as the best value budget. In contrast, D-optimality recommends a minimum budget of \$3 k to identify all uncertain parameters.

\subsection{A- and D-optimality Select Different Sensors}
As expected, A- and D-optimality objectives select different sensors.
Fig.~\ref{fig:kinetics}(b) shows that D-optimality values identify the practical non-identifiability of the model within a budget of \$2.5 k, where D-optimality values fall below 10$^{-3}$. In contrast, A-optimality values fail to show this non-estimability. The differences between A- and D-optimality can be explained through the eigendecomposition of $\bigM$, where $\lambda_j$ for $j \in {1,2,...,P}$ are the eigenvalues: 

\begin{equation}
    \text{A-optimality: trace}(\bigM) = \sum_{j=1}^{P}\lambda_j
\label{eq:trace}
\end{equation}

\begin{equation}
    \text{D-optimality: det}(\bigM) = \prod_{j=1}^{P}\lambda_j
\label{eq:det}
\end{equation}

\noindent D-optimality is the multiplication of the eigenvalues of $\bigM$, where one eigenvalue close to 0 results in an extremely small D-optimality. On the other hand, A-optimality is the summation of the eigenvalues, where an eigenvalue close to 0 has minimal impact on the A-optimality value.
The mathematical model of the system in Table \ref{tab:reactor} dictates that $C$ for one single component, even when available at different time points, cannot provide information to distinguish between the other two components and, thus, cannot identify all four parameters within Arrhenius equations. For parameter identifiability, $C$ needs to be measured for at least two components. This is seen by comparing the sensor choices shown in
Table~\ref{tab:kinetics_solution} (A-optimality) and Table~\ref{tab:kinetics_solution_D} (D-optimality).
With a budget of \$2.2 k, D-optimal solutions choose four time points from both the DCM $C_A$ and the DCM $C_B$, while A-optimality prioritizes only $C_B$ as an SCM with eight time points. Consequently, the D-optimality solution has a D-optimality value two orders of magnitude higher than the A-optimality solution at this budget.  

\subsection{Relaxed Solutions Are Not Tight}

Fig.~\ref{fig:kinetics}(b) shows the continuous relaxations are not tight, which emphasizes the need to solve the mixed-integer problems. While D-optimality MINLP solutions highlight non-estimability at low budgets, D-optimality NLP solutions incorrectly suggest the model parameters are estimable at a lowest budget of \$1 k. This is because this relaxed problem selects fractions of measurements. When the budget is below \$1.8 k, a maximum of four DCM time points can be chosen, and this limited information is insufficient to identify all parameters. SI Figure \ref{fig:kinetics_solution} futher compares the relaxed and integer optimization results.

\subsection{MO Solutions Support Heuristics and Considers Practical Constraints}
Generally, the results support the heuristics that measurements should be selected in a prioritized order from the most to the least valuable. For low budgets, the most valuable DCMs are prioritized. 
When the budget surpasses \$2 k, allowing the inclusion of SCM, the most informative SCM is chosen. 
Additional budgets are allocated to DCMs after selecting SCMs, maintaining the order of preference based on their value. This suggests that enumeration is possible when the decision space is small \citep{befort2023data}.

At a budget of \$1.4 k, the enumeration identifies that the time point combination of 7.5 and 15 minutes provides a higher D-optimality value. However, from the D-optimal solutions shown in Table \ref{tab:kinetics_solution_D}, MO selects time 7.5 and 22.5 minutes of the DCM $C_B$. 
This is because MO also considers the manual sampling limitations to ensure that, at most, one sample can be taken from the batch every ten minutes. Since time points for DCMs are 7.5 minutes apart, there is at most one time point that can be measured within each ten-minute interval. 
MO also considers various constraints, including the additional fixed costs for DCMs, the maximum sampling times for one component, and the maximum sampling times overall in the experiment. Furthermore, considering time points from DCMs usually results in high design dimensions, making exhaustive enumeration computationally expensive. Moreover, Section \ref{sec:sampling_time_constraint} in the SI demonstrates the flexibility of the framework to conduct sensitivity analysis, such as quantifying the information loss from increases the time between samples.

\subsection{Computational Aspects}

Fig.~\ref{fig:kinetics-calc} shows the computation costs of the four optimization strategies. 
A-optimality MILPs and their relaxed LPs are comparatively easier to solve, taking approximately 0.1 seconds.
On the other hand, D-optimality MINLP and its relaxed NLP problems are more challenging to solve because of their complex objective functions and possible ill-conditioning problems. At low budgets, the objective function is close to zero and $\bigM$ is (near) rank deficient, which leads to more iterations and longer solver time. Next, Fig.~\ref{fig:kinetics-calc}(a) and (c) show that with the grey-box module, the NLP problems are easily solved, typically in in 1 second with 28 \texttt{CyIpopt} iterations. Lower budgets require up to 40 \texttt{CyIpopt} iterations. 
Fig.~\ref{fig:kinetics-calc}(a) and (c) also reveal that the MINLP problems take an average of 5 seconds and 3 \texttt{MindtPy} iterations to solve, with up to 50 seconds and 35 \texttt{MindtPy} iterations at low budgets is small.  
Overall, the A-optimality MILP and LP problems solved quickly, requiring about 2\% of the solving time of D-optimality MINLP problems. While D-optimality MINLP problems required fewer iterations to solve under high budgets compared to their relaxed problems, they required approximately five-fold longer time to solve that the relaxed NLP. This is because each \texttt{MindtPy} iteration involved solving a MILP subproblem and a fixed NLP subproblem. Thus, to summarize, in this case study, A-optimality problems are easier to solve than D-optimality problems, and D-optimality relaxed problems are easier to solve than D-optimality MINLP problems.

This case study was also solved with \texttt{CVXPY} for comparison. The A-optimality MILP problem in \texttt{CVXPY} has 732 variables and 2962 constraints; it took around 2 seconds to compile and 2 seconds for \texttt{Mosek} to solve the problem. The D-optimality NLP problem had 733 variables and 2963 constraints, \texttt{CVXPY} took around 2 seconds to compile it, and \texttt{Mosek} took around 2 seconds to solve it. Both \texttt{CVXPY} and \texttt{Pyomo} generated the same results; the difference in the selection variables $\mathbf{x}$ was less than $10^{-8}$, and the difference in the computed optimality criteria was less than $10^{-3}$. Figs. \ref{fig:cvxpy_A} and \ref{fig:cvxpy_D_nlp} in Section \ref{sec:cvxpy} of the SI report the \texttt{CVXPY} results.

\section{Case Study: Rotary-Bed CO$_2$ Adsorption and Desorption}

Next, we demonstrating scalability of MO by considering a CO$_2$ adsorption and desorption system with 14 time-varying measurements which are measurable at hundreds of discretized time steps.

\subsection{Problem Statement} 

\begin{figure}[!ht]
\centering
\includegraphics[width=0.95\textwidth]{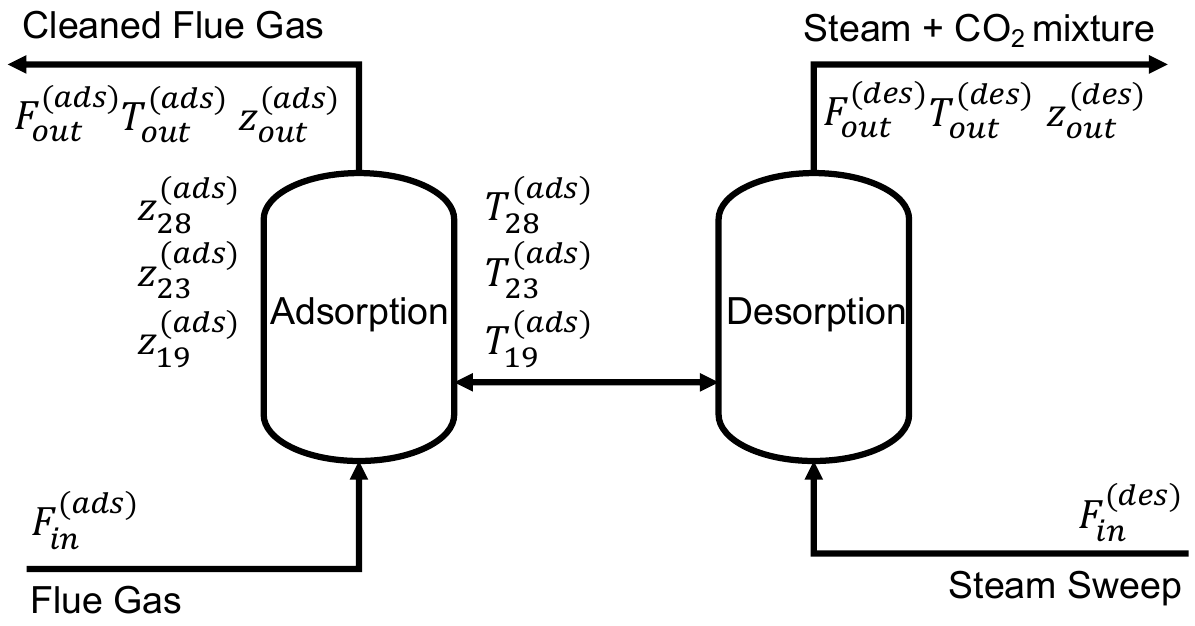}
\caption{Rotary bed adsorption-desorption system and the measurements considered for the MO problem. All 14 measurements are indexed by time (not shown). Measurements labeled with numbers (19, 23, 28) correspond to scaled positions on the column.} 
\label{fig:rotary-set}
\end{figure}

The CO$_2$ adsorption and desorption system illustrated in Fig.~\ref{fig:rotary-set} uses a functionalized metal-organic framework (MOF) to separate CO$_2$ via adsorption \citep{ezeobinwune2020modeling}. The reactive rotary packed bed (RPB) is a rotating cylindrical wheel divided into two equal sections: adsorption and regeneration (desorption of CO$_2$). The solid matrix is attached to a central rotor, which rotates at a constant speed. In this RPB, CO$_2$ from the flue gas reacts with the impregnated diamine in the adsorption section, while CO$_2$ is released by temperature swing in the regeneration section. The regeneration and adsorption sections are in a co-current flow arrangement with respect to their streams. 
To achieve a large swing in the CO$_2$ loading of the sorbent, a static embedded cooler is used in the capture section, and a static embedded heater is used in the desorption section where steam is used as the heating utility. The process model, which includes coupled mass and heat transfer with reaction kinetics, was initially developed in Aspen Custom Modeler (version 11) by \citet{ezeobinwune2020modeling}. We consider a reactor with a diameter of 3m and length of 3m, which accommodates an inlet flue gas flow rate of 365 kmol$\cdot$hr$^{-1}$.
 
 The adsorption bed is discretized in the axial direction to have 33 positions, where position 0 is the inlet and position 32 is the outlet. Finer discretizations are used near the inlet of the bed to capture the significant profile changes caused by rapid kinetics near the inlet. Between a normalized length of 0 and 0.1, there are 17 positions distributed evenly. The rest of the bed, from 0.1 to 1, contains 16 nodes distributed evenly. In addition, the model is also discretized in the angular direction with 17 total nodes. The angular distance is normalized by a half turn through the adsorption section, where position 0 is the solids inlet, and 0.5 is half a rotation through the adsorption section. Between a normalized length of 0 and 0.05, there are seven nodes distributed evenly. The rest of the bed, from 0.05 to 0.5, has ten nodes distributed evenly. Finer discretizations are used to capture the fast kinetics changes in the solids input side and thus help with numerical convergence of the model.
 
The goal of MO is to decide the measurement strategy to maximize the precision of the model parameters --- the mass transfer coefficient $MTC$, the heat transfer coefficient $HTC$, and heat adsorption parameters $DH$, $Iso_1$, and $Iso_2$ --- from 14 candidate measurements: flowrates at gas inlets $F_{in}^{(ads)}$ and  $F_{in}^{(des)}$, and gas outlets $F_{out}^{(ads)}$ and $F_{out}^{(des)}$; temperatures at the gas outlets $T_{out}^{(ads)}$ and $T_{out}^{(des)}$, and the axial positions 19, 23, and 28 inside the adsorption bed $T_{19}^{(ads)}$, $T_{23}^{(ads)}$, $T_{28}^{(ads)}$; gas compositions at the gas outlets $z_{out}^{(ads)}$ and  $z_{out}^{(des)}$, and the axial positions 19, 23, and 28 inside the adsorption bed $z_{19}^{(ads)}$, $z_{23}^{(ads)}$, $z_{28}^{(ads)}$. (These axial positions correspond to 0.80 m, 1.48 m, and 2.33 m from the inlet.) All the time-varying measurements may be sampled at 110 time points distributed evenly in 220 minutes. The observation errors are assumed to be independent both in time and across measurements, with a constant standard deviation of 1 K for $T$, 1 mol$^{-1}$ L for $F$, and 0.01 for $z$. The dynamic sensitivities are computed via central finite difference with 1\% perturbation of the model parameters.

\begin{table}[htb!]
\centering
\caption{Measurement costs and types for the rotary bed case study.} 
\begin{tabular}{| l c c |}
\toprule
Name & Installation (\$) & Measurement cost (\$/sample) \\
\midrule
$F_{in}^{(ads)}$ & 1000 & 0 \\ 
$F_{out}^{(ads)}$ & 1000 & 0 \\ 
$T_{out}^{(ads)}$ & 500 & 0 \\ 
$F_{in}^{(des)}$ & 1000 & 0 \\ 
$F_{out}^{(des)}$ & 1000 & 0 \\ 
$T_{out}^{(des)}$ & 500 & 0 \\ 
$T_{19}^{(ads)}$ & 1000 & 0 \\ 
$T_{23}^{(ads)}$ & 1000 & 0 \\ 
$T_{28}^{(ads)}$ & 1000 & 0 \\ 
$z_{out}^{(ads), SCM}$ & 7000 & 0 \\  
$z_{out}^{(des), SCM}$ & 7000 & 0 \\ 
$z_{out}^{(ads), DCM}$ & 100 & 100 \\  
$z_{out}^{(des), DCM}$ & 100 & 100 \\
$z_{19}^{(ads)}$ & 500 & 100 \\ 
$z_{23}^{(ads)}$ & 500 & 100 \\ 
$z_{28}^{(ads)}$ & 500 & 100 \\
\bottomrule
\end{tabular}
\label{tab:rotary_price}
\end{table}

Table \ref{tab:rotary_price} shows the installation and measurement costs for this experimental system. According to \citet{liptak2003instrument}, thermal mass flowmeters are an economical choice for the four flowrates, costing around \$1 k per sensor. Thermocouples are used for temperature measurements. We consider measuring $z_{out}^{(ads)}$ and $z_{out}^{(des)}$ either with online GC machines as SCMs, or with manual sampling as DCMs.

Furthermore, considering human resources constraints, a maximum of five time points can be chosen for each DCM, and a maximum of 20 time points in total are allowed for all DCMs. Additionally, the selecting time points must maintain a minimum separation of at least ten minutes between each other.

Four optimization strategies are solved: maximizing A-optimality of the MILP problem, its relaxed LP problem; maximizing D-optimality of the MINLP problem, and its relaxed NLP problem.
The A-optimality LP problem has 157,663 continuous variables, 17 equality constraints, and 471,914 inequality constraints. 
The A-optimality MILP problem has 17 continuous variables, 157,646 binary variables, 17 equality constraints, and 471,914 inequality constraints. 
The D-optimality NLP problem has 157,679 variables, 33 equality constraints, and 471,914 inequality constraints. 
The D-optimality MILP problem has 33 continuous variables, 157,646 binary variables, 33 equality constraints, and 471,914 inequality constraints. 

\setlength\tabcolsep{3pt} 
\begin{longtable}{|l | c | c | c | c | c | c| c| c |c |c |c |c |c |c |c |c |c |c |c |c |c |c |c |c |c| }
\caption{Optimal solutions of SCMs for the A-optimality MILPs where each column corresponds to a measurement budget of \$1 k to \$25 k. '1' and '0' indicate if each SCM is selected or not selected.}\\
\hline
\multicolumn{1}{|c|}{Budget [\$ k]} & 
\multicolumn{1}{c|}{\textbf{1}} & \multicolumn{1}{c|}{\textbf{2}} & \multicolumn{1}{c|}{\textbf{3}} &
\multicolumn{1}{c|}{\textbf{4}} &
\multicolumn{1}{c|}{\textbf{5}} &
\multicolumn{1}{c|}{\textbf{6}} &
\multicolumn{1}{c|}{\textbf{7}} &
\multicolumn{1}{c|}{\textbf{8}} &
\multicolumn{1}{c|}{\textbf{9}} &
\multicolumn{1}{c|}{\textbf{10}} &
\multicolumn{1}{c|}{\textbf{11}} &
\multicolumn{1}{c|}{\textbf{12}}&
\multicolumn{1}{c|}{\textbf{13}}&
\multicolumn{1}{c|}{\textbf{14}}&
\multicolumn{1}{c|}{\textbf{15}}&
\multicolumn{1}{c|}{\textbf{16}}&
\multicolumn{1}{c|}{\textbf{17}}&
\multicolumn{1}{c|}{\textbf{18}}&
\multicolumn{1}{c|}{\textbf{19}}&
\multicolumn{1}{c|}{\textbf{20}}&
\multicolumn{1}{c|}{\textbf{21}}&
\multicolumn{1}{c|}{\textbf{22}}&
\multicolumn{1}{c|}{\textbf{23}}&
\multicolumn{1}{c|}{\textbf{24}}&
\multicolumn{1}{c|}{\textbf{25}} \\ \hline 
\endfirsthead

\multicolumn{26}{c}%
{{\bfseries \tablename\ \thetable{} -- continued from previous page}} \\
\hline 
\multicolumn{1}{|c|}{Budget [\$ k]} & 
\multicolumn{1}{c|}{\textbf{1}} & \multicolumn{1}{c|}{\textbf{2}} & \multicolumn{1}{c|}{\textbf{3}} &
\multicolumn{1}{c|}{\textbf{4}} &
\multicolumn{1}{c|}{\textbf{5}} &
\multicolumn{1}{c|}{\textbf{6}} &
\multicolumn{1}{c|}{\textbf{7}} &
\multicolumn{1}{c|}{\textbf{8}} &
\multicolumn{1}{c|}{\textbf{9}} &
\multicolumn{1}{c|}{\textbf{10}} &
\multicolumn{1}{c|}{\textbf{11}} &
\multicolumn{1}{c|}{\textbf{12}}&
\multicolumn{1}{c|}{\textbf{13}}&
\multicolumn{1}{c|}{\textbf{14}}&
\multicolumn{1}{c|}{\textbf{15}}&
\multicolumn{1}{c|}{\textbf{16}}&
\multicolumn{1}{c|}{\textbf{17}}&
\multicolumn{1}{c|}{\textbf{18}}&
\multicolumn{1}{c|}{\textbf{19}}&
\multicolumn{1}{c|}{\textbf{20}}&
\multicolumn{1}{c|}{\textbf{21}}&
\multicolumn{1}{c|}{\textbf{22}}&
\multicolumn{1}{c|}{\textbf{23}}&
\multicolumn{1}{c|}{\textbf{24}}&
\multicolumn{1}{c|}{\textbf{25}} \\ \hline 
\endhead

\hline \multicolumn{26}{|r|}{{Continued on next page}} \\ \hline
\endfoot

\hline \hline
\endlastfoot
 
$F_{in}^{ads}$ & 0  & 0  & 0  & 1 & 1& 1& 1
& 1& 1& 1& 1& 1& 1& 1& 1& 1& 1& 1& 1& 1& 1 & 1& 1& 1& 1 \\ 
$F_{out}^{ads}$ & 0  & 0  & 0  & 0 & 0& 0& 0
& 0& 0& 1& 1& 0& 0& 0& 0& 0& 1& 1& 1& 0& 0 & 0& 1& 1& 1 \\
$T_{out}^{ads}$ & 1  & 1  & 1  & 1 & 1& 1& 1
& 1& 1& 1& 1& 1& 1& 1& 1& 1& 1& 1& 1& 1& 1 
& 1& 1& 1& 1 \\
$F_{in}^{des}$ & 0  & 0  & 1  & 1 & 1& 1& 1
& 1& 1& 1& 1& 1& 1& 1& 1& 1& 1& 1& 1& 1& 1 
& 1& 1& 1& 1 \\
$F_{out}^{des}$ & 0  &1  & 1  & 1 & 1& 1& 1
& 1& 1& 1& 1& 1& 1& 1& 1& 1& 1& 1& 1& 1& 1 
& 1& 1& 1& 1 \\
$T_{out}^{des}$ & 1  & 1  & 1  & 1 & 1& 1& 1
& 1& 1& 1& 1& 1& 1& 1& 1& 1& 1& 1& 1& 1& 1 
& 1& 1& 1& 1 \\
$T_{19}^{ads}$ & 0  & 0  & 0  & 0 & 1& 1& 1
& 1& 1& 1& 1& 1& 1& 1& 1& 1& 1& 1& 1& 1& 1 
& 1& 1& 1& 1 \\
$T_{23}^{ads}$ & 0  & 0  & 0  & 0 & 0& 1& 1
& 1& 1& 1& 1& 0& 1& 1& 1& 1& 1& 1& 1& 1& 1 
& 1& 1& 1& 1 \\
$T_{28}^{ads}$ & 0  & 0  & 0  & 0 & 0& 0& 0
& 1& 1& 1& 1& 0& 0& 0& 1& 1& 1& 1& 1& 0& 0 
& 1& 1& 1& 1 \\
$z_{out}^{ads}$ & 0  & 0  & 0  & 0 & 0& 0& 0
& 0& 0& 0& 0& 1& 1& 1& 1& 1& 1& 1& 1& 1& 1 
& 1& 1& 1& 1 \\
$z_{out}^{des}$ & 0  & 0  & 0  & 0 & 0& 0& 0
& 0& 0& 0& 0& 0& 0& 0& 0& 0& 0& 0& 0& 1& 1 
& 1& 1& 1& 1 
\label{tab:rotary_solution_SCM_A}
\end{longtable}

\setlength\tabcolsep{3pt} 
\begin{longtable}{|l | c | c | c | c | c | c| c| c |c |c |c |c |c |c |c |c |c |c |c |c |c |c |c |c |c| }
\caption{Optimal solutions of DCMs for the A-optimality MILPs where each column corresponds to a measurement budget of \$1k to \$25k. The selected sample times [min] of DCMs are reported.} \\
\hline
\multicolumn{1}{|c|}{Budget [\$k]} & 
\multicolumn{1}{c|}{\textbf{7}} &
\multicolumn{1}{c|}{\textbf{8}} &
\multicolumn{1}{c|}{\textbf{9}} &
\multicolumn{1}{c|}{\textbf{10}} &
\multicolumn{1}{c|}{\textbf{11}} &
\multicolumn{1}{c|}{\textbf{14}}&
\multicolumn{1}{c|}{\textbf{15}}&
\multicolumn{1}{c|}{\textbf{16}}&
\multicolumn{1}{c|}{\textbf{17}}&
\multicolumn{1}{c|}{\textbf{18}}&
\multicolumn{1}{c|}{\textbf{19}}&
\multicolumn{1}{c|}{\textbf{21}}&
\multicolumn{1}{c|}{\textbf{22}}&
\multicolumn{1}{c|}{\textbf{23}}&
\multicolumn{1}{c|}{\textbf{24}}&
\multicolumn{1}{c|}{\textbf{25}}\\ \hline 
\endfirsthead

\multicolumn{26}{c}%
{{\bfseries \tablename\ \thetable{} -- continued from previous page}} \\
\hline 
\multicolumn{1}{|c|}{Budget [\$k]} & 
\multicolumn{1}{c|}{\textbf{7}} &
\multicolumn{1}{c|}{\textbf{8}} &
\multicolumn{1}{c|}{\textbf{9}} &
\multicolumn{1}{c|}{\textbf{10}} &
\multicolumn{1}{c|}{\textbf{11}} &
\multicolumn{1}{c|}{\textbf{14}}&
\multicolumn{1}{c|}{\textbf{15}}&
\multicolumn{1}{c|}{\textbf{16}}&
\multicolumn{1}{c|}{\textbf{17}}&
\multicolumn{1}{c|}{\textbf{18}}&
\multicolumn{1}{c|}{\textbf{19}}&
\multicolumn{1}{c|}{\textbf{21}}&
\multicolumn{1}{c|}{\textbf{22}}&
\multicolumn{1}{c|}{\textbf{23}}&
\multicolumn{1}{c|}{\textbf{24}}&
\multicolumn{1}{c|}{\textbf{25}}\\ \hline 
\endhead

\hline \multicolumn{26}{|r|}{{Continued on next page}} \\ \hline
\endfoot

\hline \hline
\endlastfoot
 
$z_{out}^{ads}$ &\ 
&\ &126 &126 &124  &\ 
&\ &\ &\ &\ &\ &\
&\ &\ &\ &\ 
\\
 &\ 
&\ &136 &136 &134  &\ 
&\ &\ &\ &\ &\  &\
&\ &\ &\ &\ 
\\
 &\
&\ &180 &182 &180  &\ 
&\ &\ &\ &\ &\  &\
&\ &\ &\ &\ 
\\
 &\
&\ &190 &192 &190  &\ 
&\ &\ &\ &\ &\  &\
&\ &\ &\ &\ 
\\
&\ 
&\ &\ &\ &200 &\ 
&\ &\ &\ &\ &\  &\
&\ &\ &\ &\ 
\\
$z_{out}^{des}$   &\ 
 &\ &2 &2 &2 &\ 
 &\ &2 &2 &2 &2  &\
 &\ &\ &\ &\ \\

 &\ 
 &\ &12 &12 &12 &\ 
   &\ &12 &12 &12 &12  &\ 
   &\ &\ &\ &\ \\

 &\ 
 &\ &22 &22 &22  &\
  &\ &22 &22 &22 &22  &\
  &\ &\ &\ &\
  \\

&\ 
 &\ &32 &32 &32  &\ 
  &\ &32 &32 &32 &32  &\
  &\ &\ &\ &\
   \\

 &\ 
 &\ &42 &42 &42  &\ 
  &\ &42 &42 &42 &42  &\
  &\ &\ &\ &\
  \\

$z_{19}^{ads}$   &\
&\ &\ &\ &\ &\ 
 &\ &\ &\ &158 &76  &\ 
  &\ &\ & 168 &76 \\

 &\
 &\ &\ &\ &\ &\ 
 &\ &\ &\ &186 &140  &\ 
  &\ &\ & 180 &140 \\

 &\
  &\ &\ &\ &\ &\ 
  &\ &\ &\ &196 &150  &\ 
  &\ &\ & 196 &150\\

 &\
   &\ &\ &\ &\ &\ 
   &\ &\ &\ &206 &180  &\ 
   &\ &\ & 206 &180\\

 &\
    &\ &\ &\ &\ &\  
    &\ &\ &\ &216 &190  &\ 
    &\ &\ & 216 &190\\

$z_{23}^{ads}$   &\ 
 &\ &\ &\ &158  &\ 
  &\ &\ &\ &\ &160  &\ 
   &\ &\ &\ &160\\

&\ 
  &\ &\ &\ &210 &\ 
   &\ &\ &\ &\ &170  &\ 
    &\ &\ &\ &170 \\

  &\ 
   &\ &\ &\ &220 &\ 
    &\ &\ &\ &\ &200  &\ 
     &\ &\ &\ &200 \\

&\ 
    &\ &\ &\ &\  &\ 
     &\ &\ &\ &\ &210  &\ 
      &\ &\ &\ &210 \\

&\ 
     &\ &\ &\ &\  &\ 
      &\ &\ &\ &\ &220  &\ 
       &\ &\ &\ &220 \\
    
$z_{28}^{ads}$  &86 
&86 &86 &86 &84 &86 
&86 &86 &86 &86 &86  &86
&86 &86 &86 &86\\

 &96 
 &96 &96 &96 &94  &96
 &96 &96 &96 &96 &96  &96
 &96 &96 &96 &96\\

 &106 
  &106 &106 &106 &104  &106
  &106 &106 &106 &106 &106  &106
   &106 &106 &106 &106\\

  &116 
   &116 &116 &116 &114  &116
   &116 &116 &116 &116 &116  &116
   &116 &116 &116 &116\\

   &126 
    &126 &\ &\ &170  &126 
    &126 &126 &126 &126 &126 &126
    &126 &126 &126 &126 
\label{tab:rotary_solution_DCM_A}
\end{longtable}

\begin{figure}[!ht]
\centering
\includegraphics[width=0.95\textwidth]{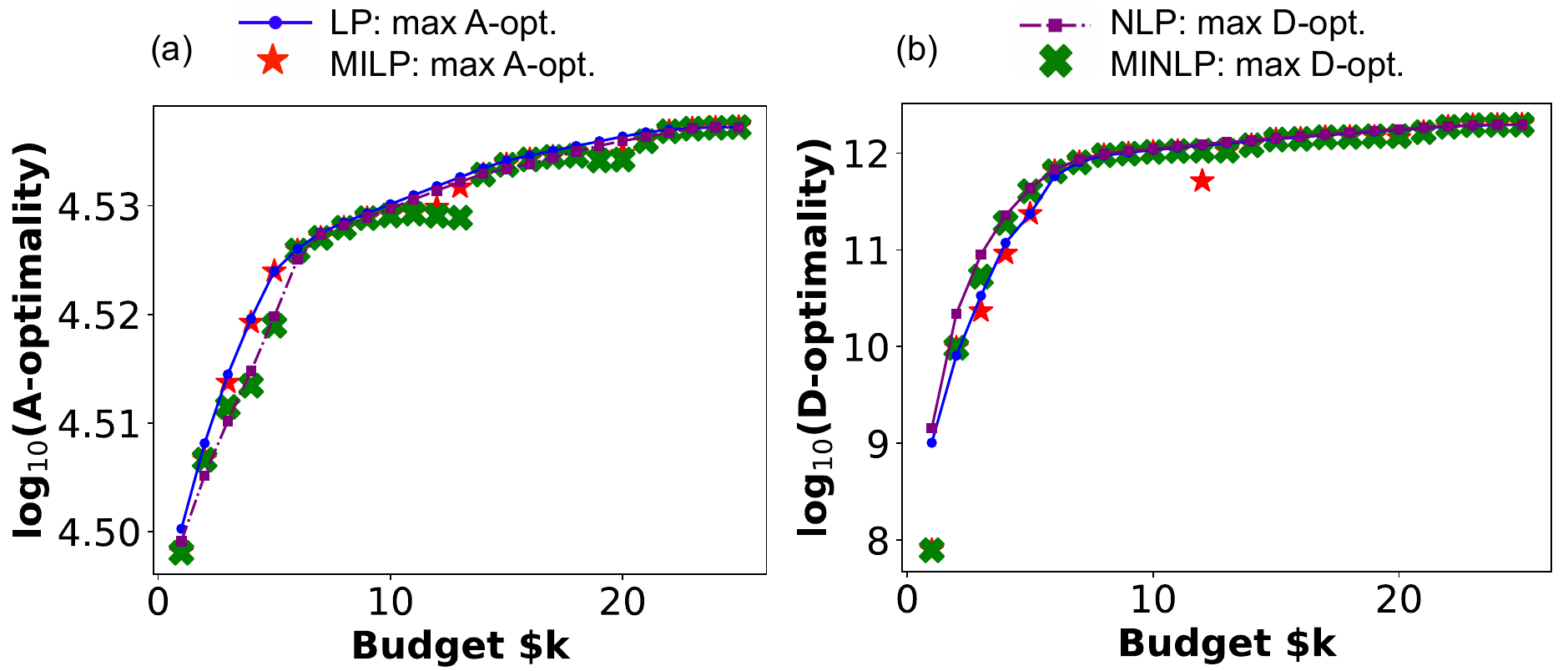}
\caption{Pareto-optimal trade-off between measurement budgets versus (a) A-optimality (trace of FIM) and (b) D-optimality (determinant of FIM) for the CO$_2$ capture case study considering four optimization strategies: maximizing A-optimality of the MILP problem (red stars), its relaxed LP problem (blue line); maximizing D-optimality of the MINLP problem (green crosses), and its relaxed NLP problem (purple line). In (a), the blue line is an upper bound for the red stars, while in (b) the purple line is an upper bound for the green crosses.} 
\label{fig:rotary}
\end{figure}

\begin{figure}[!ht]
\centering
\includegraphics[width=\textwidth]{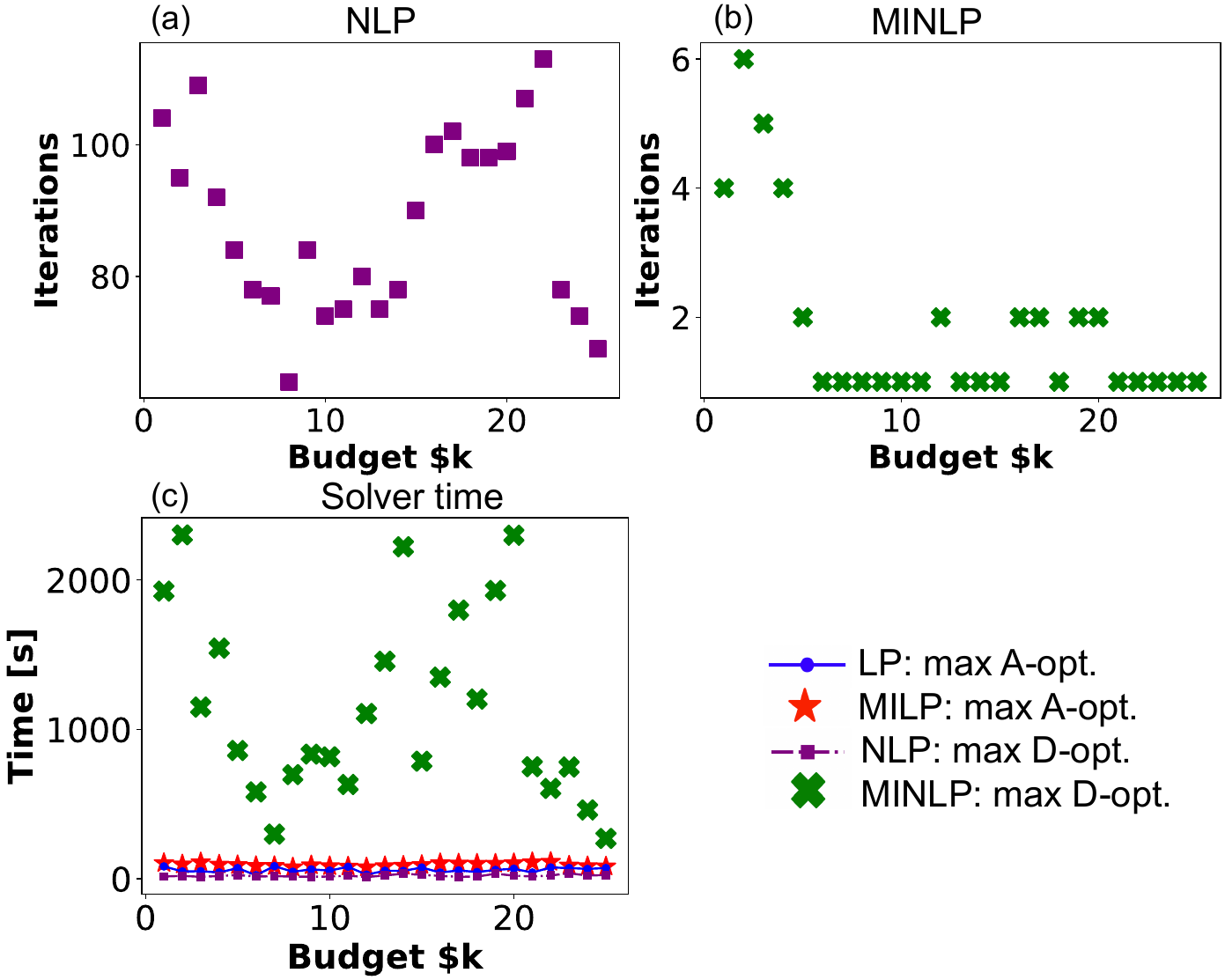}
\caption{Computational results for the CO$_2$ capture case study at different budgets. (a) the number of \texttt{CyIpopt} iterations of the maximize D-optimality NLPs, (b) the number of \texttt{MindtPy} iterations for the maximize D-optimality MINLPs, and (c) the computational time for all four optimization strategies. The optimization problems are solved from the lowest to the highest budget.} 
\label{fig:rotary-calc}
\end{figure}

\subsection{Pareto-Optimal Trade-Offs Between Measurement Budgets and Information}

Fig.~\ref{fig:rotary} assists practitioners in choosing the most appropriate measurement budget by visualizing the Pareto-optimal trade-offs between measurement costs and the A-optimality and D-optimality metrics. Fig.~\ref{fig:rotary} was generated by solving four optimization problems over the budget range of \$1k to \$25k. 
The A-optimality MILP solutions are shown in Table \ref{tab:rotary_solution_SCM_A} and Table \ref{tab:rotary_solution_DCM_A}; the D-optimality MINLP solutions are reported in Table~\ref{sec:rotaryD} of the SI.
Fig.~\ref{fig:rotary}(a) shows the A-optimality LP solutions are the upper bounds of A-optimality MILP solutions. For comparison, Fig.~\ref{fig:rotary}(a) also shows the A-optimality metrics for the D-optimality NLP and MINLP solutions. As expected, the D-optimality NLP does not provide an upper bound for the A-optimality metric for the D-optimality MINLP in Fig.~\ref{fig:rotary}(a) because these two optimization problems consider a different objective. Similar observations also apply to Fig.~\ref{fig:rotary}(b) for evaluating the D-optimality metric for the A-optimal solutions.
All four solutions show a sharp increase in A- and D-optimality from \$1k to \$5k, followed by a much more gradual increase from \$6k to \$25k. A practitioner may conclude from Fig.~\ref{fig:rotary} that \$5k is a recommended value budget. Even with the lowest budget \$1k, D-optimality is 10$^{8}$, indicating that all parameters in this problem are estimable.

\subsection{A- and D-optimality Select Different Sensors}

The A- and D-optimality objectives are more closely related in this case study compared to the reaction kinetics one. Specifically, Fig.~\ref{fig:rotary}(a) shows that maximizing D-optimality also improves A-optimality values, and Fig.~\ref{fig:rotary}(b) shows that maximizing A-optimality also improves D-optimality values. 

However, the A-optimality MILP results show a low D-optimality value at a budget of \$12k, nearly 0.5 orders of magnitude smaller than the D-optimality solution. With this budget, the A-optimality solution includes an expensive SCM $z_{out}^{(ads)}$, while the D-optimality solution includes other cheaper SCMs and DCMs, and drops $z_{out}^{(ads)}$, which leads to a higher D-optimality value. This indicates that A-optimality cannot consistently perform well from a D-optimality perspective. This is because A-optimality relies solely on the diagonal elements of $\bigM$, whereas D-optimality considers correlations between parameters.

\subsection{Relaxed Solutions Are Not Tight}

In this case study, the relaxed problems find integer solutions for SCMs but select fractions of time points of DCMs in solutions. For example, for the D-optimality problem at a budget of \$6k, instead of choosing one time point from $z_{28}^{(ads)}$, the NLP relaxation chooses five time points, each with a fraction of 0.2. This shows that the relaxed MO can only provide a preliminary analysis of this problem and cannot always lead to the same observations as the mixed-integer problems. This highlights that the mixed-integer problem is still crucial for choosing sensors in this system.

\subsection{MO Solutions Support Heuristics and Consider Practical Constraints}
Table \ref{tab:rotary_solution_SCM_A} shows the A-optimality MILP solutions under different budgets. The solutions support the heuristics that measurements should be selected in a prioritized order from the most to the least valuable. For low budgets (less than \$6k in this case), SCMs that cost comparatively less but contain substantial information are prioritized. 
Within the intermediate budgets ranging from \$7k to \$20k, DCMs are gradually added to  utilize the available budget.
When the budgets continue to increase, the solutions tend to select expensive SCMs, sometimes at the expense of deselecting certain DCMs and SCMs. For instance, switching the budget from \$11k to \$12k, MO removes the selection of three SCMs and all time points of DCMs to allow the selection of the SCM $z_{out}^{(ads)}$. SI Section \ref{sec:rotaryA} further elaborates on the physical interpretation of these results.

However, practical constraints and experimental requirements make enumeration impractical and require optimization frameworks to solve this problem. For example, due to manual sampling limitations, a maximum of five time points can be chosen from each DCM. Table \ref{tab:rotary_solution_DCM_A} shows that this constraint is satisfied by the optimal solutions at all budgets. For instance, at a budget of \$16k, while selecting more time points from $z^{ads}_{28}$ can provide a higher D-optimality value, MO selects time points from $z^{des}_{out}$ since the number of chosen time points from $z^{ads}_{28}$ has reached the limit of five. 
MO also considers various constraints, including the additional fixed costs for DCMs, the maximum sampling times for one component, and the maximum sampling times overall in the experiment. In this case, the number of candidate sensors and time points reaches 561, making exhaustive enumeration computationally expensive. This means that optimization is crucial to solving this problem under practical requirements. 

\subsection{Computational Aspects}

Finally, Fig.~\ref{fig:rotary-calc} shows the computational expense of the four optimization strategies. The A-optimality relaxation takes 57.8 seconds on average to solve. Similarly, the A-optimality MILPs take 19.7 seconds on average, which takes less time than the LPs because the MILPs are initialized using the solution of the LP relaxations. In contrast, the D-optimality MINLPs and their relaxed NLPs are more challenging to solve with the grey-box module. Fig.~\ref{fig:rotary-calc}(a) shows that NLP problems are solved on average in 93 seconds with 88 \texttt{CyIpopt} iterations. The variations in the number of iterations required to solve the problem at different budgets are influenced by several factors. Selecting SCMs becomes relatively easier with higher budgets since there are more SCMs available to choose from a fixed total number of SCMs, resulting in fewer iterations. On the other hand, the challenge of selecting DCMs may increase since, as the number of time points increases, the measurement space expands, leading to more iterations. This fluctuation is also caused by the initialization strategy we adopted. We applied the A-optimality MILP solutions as the initial points for NLP problems, while this initial point is not close to the optimal solution since NLP solutions usually incorporate fractions of measurements, while MILP solutions only have binary decisions. Therefore, the number of NLP iterations is non-monotonic since the proximity of the initial point to the optimal solution can be comparatively different for different budgets.
Fig.~\ref{fig:rotary-calc}(b) shows that the MINLP problems take an average of 1144 seconds and 2 \texttt{MindtPy} iterations to solve. Leveraging the good starting points provided by A-optimality MILP solutions, MINLP problems take relatively few iterations to converge. Interestingly, the number of iterations needed for solving an MINLP problem decreases with an increasing budget. This trend can be explained by the fact that at low budgets, there is more flexibility in choosing SCMs from a fixed total number of SCMs, leading to larger number of iterations. Moreover, as the budget increases, the MINLP solutions are more similar to the initial points of the MILP solution, leading to less iterations. 
In each iteration, the solver time for the master MILP problem is approximately 200 seconds, and the solver time for the fixed NLP problem is around 0.1 seconds. Given the large problem size from the many integer variables and inequality constraints, the wall clock solving time for each iteration is around 400 seconds, where the master problem takes roughly 300 seconds to formulate and solve, and the fixed NLP problem takes around 100 seconds to formulate and solve. 

In summary, this case study demonstrates the scalability of the MO model and confirms general trends about the computational difficult of the four optimization problems. The A-optimality MILP and LP problems generally solve much faster than D-optimality problems, requiring less than half the solving time of D-optimality NLP problems and only 4\% of the solving time of D-optimality MINLP problems. D-optimality MINLP problems take around 12 times longer to solve than their relaxed NLP problems. For this large-scale case study, A-optimality problems are easier to solve than D-optimality problems, and D-optimality relaxed problems are easier to solve than D-optimality MINLP problems.

\section{Conclusions}

The selection of measurements is critical for model-based data analytics and to construct predictive digital twins; decisions about instrumentation and equipment design are often made months to years before experimental campaigns begin. Measurement optimization provides science and engineering communities with a systematic approach for selecting the most informative measurements within certain budgets and requirements. This capability accelerates collaborations (e.g., it was inspired by questions arising from \cite{ouimet2022data}) and promotes the development of predictive digital twins. 

This paper presents a convex MINLP formulation to calculate the Pareto-optimal trade-off between information content, measured by A- and D-optimality, and measurement budgets. The optimization framework is highly flexible and easily accommodates practical constraints such as a minimum time between manual measurements. The key computational contributions of this work are evaluating the D-optimality objective function and its gradient using \texttt{SciPy}, which is incorporated into \texttt{Pyomo} using the new \texttt{ExternalGreyBoxModel} module. A modified version of the MINLP solver \texttt{MindtPy} is developed and demonstrated for use with grey-box constraints. The D-optimality MINLP problem of the large-scale rotary bed system, with 157,646 binary variables and 471,914 inequality constraints, can be solved in an average of 2 \texttt{MindtPy} iterations and 1144 seconds total, showing that scalability of the proposed computational approach. The problem is also formulated and solved using \texttt{CVXPY} to verify the optimal solutions of the convex formulations. 

Several general observations can be drawn from the reaction kinetics and CO$_2$ capture case studies. In both applications, A- and D-optimality choose different sensors driven by different objectives. Notably, D-optimality recognizes practical identifiability issues with small budgets and selects measurements to avoid this, whereas A-optimality does not. Furthermore, many of the relaxed problem solutions were not tight and consequently lost their ability to detect estimability issues, highlighting the importance of solving the MI(N)LP, especially with limited budgets. Both case studies support the heuristics to prioritize measurements from most to least valuable unless constraints prevent this. From a computational perspective, A-optimality problems are easier to solve than D-optimality problems, and D-optimality relaxed problems are easier to solve than D-optimality MINLP problems.

Possible future extensions include the deployment of other MBDoE criteria for certain analysis goals. G-optimality, which minimizes the maximum variance of any predicted value \citep{wong1994comparing}, and V-optimality, which minimizes the prediction variance \citep{pronzato1989experiment}, can be used to consider the impact of parameter uncertainty on response predictions. E-optimality \citep{franceschini2008model}, which maximizes the minimum eigenvalue of FIM, can be optimized to increase the practical identifiability of parameters. Furthermore, there is potential for adapting this problem formulation to multi-objective optimization of multiple information content metrics, costs, and model discrimination goals. Using exact Hessian information of the D-optimality objective may further improve computational efficiency. Additionally, the problem can be extended to incorporate surrogate modeling techniques and Bayesian hybrid models (BHMs) \citep{eugene2023learning} to consider both parameter and model-form uncertainties. Furthermore, the proposed framework can be extended to consider the joint optimization of measurement selection and experimental conditions.
Finally, the new capability to calculate log determinants in Pyomo via the \texttt{ExternalGreyBoxModel} module widely applies to many optimization problems in data science and applied statistics with nonlinear models, e.g., maximum likelihood estimation and MBDoE.

\section*{CRediT Authorship Contribution Statement}

Jialu Wang: Conceptualization, Methodology, Software, Investigation, Formal Analysis, Data creation, Validation, Visualization, Writing - original draft.
Zedong Peng: Software.  
Ryan Hughes: Process model simulation. 
Debangsu Bhattacharyya: Supervision. 
David E. Bernal Neira: Software, Supervision.
Alexander W. Dowling: Conceptualization, Methodology, Software, Project administration, Resources, Supervision, Writing - review \& editing,  Funding acquisition.

\section*{Data Availability Statement}

The data that support the findings of this study are openly available at \url{https://github.com/dowlinglab/measurement-opt}.

\section*{Acknowledgements}

The authors graciously acknowledge funding from the U.S. Department of Energy, Office of Fossil Energy and Carbon Management, through the Carbon Capture Program. 

\section*{The Disclaimer}

This project was funded by the Department of Energy, National Energy Technology Laboratory an agency of the United States Government, through a support contract. Neither the United States Government nor any agency thereof, nor any of its employees, nor the support contractor, nor any of their employees, makes any warranty, expressor implied, or assumes any legal liability or responsibility for the accuracy, completeness, or usefulness of any information, apparatus, product, or process disclosed, or represents that its use would not infringe privately owned rights. Reference herein to any specific commercial product, process, or service by trade name, trademark, manufacturer, or otherwise does not necessarily constitute or imply its endorsement, recommendation, or favoring by the United States Government or any agency thereof. The views and opinions of authors expressed herein do not necessarily state or reflect those of the United States Government or any agency thereof.

\section*{Conflict of Interest Statement}
\noindent Nothing to declare.

\newpage
\section*{Nomenclature}

\begin{table}[htb!] 
    \small
    \caption{Latin letters (Reaction kinetics case study)}
    \label{tab:notation-case}
    \begin{tabular}{l l }
    $A_1, A_2$ & Arrhenius equation pre-exponential factors, 1/h \\
    $C_A, C_B, C_C$ & Concentration of components, mol/L \\ 
    $C_{A0}$        & Initial concentration of component A, mol/L \\
    $E_1, E_2$ & Arrhenius equation activation energies, kJ/mol \\ 
   $k_1, k_2$ & Rate constants, 1/h \\
    \end{tabular}
\end{table}

\begin{table}[htb!] 
    \small
    \caption{Latin letters (Rotary-bed case study)}
    \label{tab:notation-case2}
    \begin{tabular}{l l }
   $DH$ & Adjustment factor for heat of adsorption \\ 
   $F_{in}^{(ads)}$ & Inlet gas flowrate of the adsorption tower, mol/s \\
   $F_{in}^{(des)}$ & Inlet gas flowrate of the desorption tower, mol/s \\ 
   $F_{out}^{(ads)}$ & Outlet gas flowrate of the adsorption tower, mol/s \\ 
   $F_{out}^{(des)}$ & Outlet gas flowrate of the desorption tower, mol/s\\ 
   $HTC$ & Heat transfer coefficient for gas-to-solid heat transfer \\ 
   $Iso_1$ & Adjustment factor for the overall prediction of the isotherm model \\ 
   $Iso_2$ & Adjustment factor for a temperature-dependent parameter in the isotherm model \\
   $MTC$ & Overall mass transfer coefficient \\
   $T_{out}^{(ads)}$  & Outlet gas temperature of the adsorption tower, K \\
   $T_{out}^{(des)}$ & Outlet gas temperature of the desorption tower, K \\
   $T_{19}^{(ads)}$  & Gas temperature of the adsorption tower at axial position index 19, K \\
   $T_{23}^{(ads)}$  & Gas temperature of the adsorption tower at axial position index 23, K \\
   $T_{28}^{(ads)}$  & Gas temperature of the adsorption tower at axial position index 28, K \\
   $z_{out}^{(ads)}$  & Outlet gas fraction of the adsorption tower \\
   $z_{out}^{(des)}$  & Outlet gas fraction of the desorption tower \\
   $z_{19}^{(ads)}$  & Gas fraction of the adsorption tower at axial position index 19 \\
   $z_{23}^{(ads)}$  & Gas fraction of the adsorption tower at axial position index 23 \\
   $z_{28}^{(ads)}$ & Gas fraction of the adsorption tower at axial position index 28 \\

    \end{tabular}
\end{table}

\begin{table}[htb!] 
    \small
    \caption{Latin Letters (General method)}
    \label{tab:nom}
    \begin{tabular}{l l }
    $B$ & Budget [\$] \\
    $c$ & Fixed or per-measurement cost \\
    $C$ & Number of experimental conditions \\ 
    $D$ & Number of dynamic-cost measurements \\
    $\mathcal{D}$ & Dynamic-cost measurement space \\
   $\mathbf{f}(\cdot)$ & Differential or Algebraic equations \\ 
   $K$        & Number of state variables \\
   $L_d$      & Number of time points allowed for dynamic-cost measurement $y_d$ \\ 
   $L_{total}$ & Number of all time points allowed for all dynamic-cost measurements \\
   $\mathbf{M}$& Fisher Information Matrix \\
   $P$        & Number of parameters \\
   $\mathbf{Q}$& Dynamic sensitivity matrix \\
   $S$ & Number of static-cost measurements \\
   $\mathcal{S}$ & Static-cost measurement space \\
   $\mathbf{t}$ & Time set of measurement response variables \\
   $T_{int}$ & Minimal time interval between two measured time points \\ 
   $T_k$ & Number of time points for measurement $y_k$ \\
     $x$ & Binary decision for if including a measurement \\
     $\bigy$ & Measurement response variable vector\\
     $\mathcal{Y}$ & Measurement response variable space \\
    \end{tabular}
\end{table}

\begin{table}[htb!] 
    \small
    \caption{Greek letters}
    \label{tab:not-math}
    \begin{tabular}{l l }
        $\bigtheta$  & Parameter vector \\ 
        $\Theta$ & Parameter space \\ 
        $\bigphi$ & Design variable vector \\
        $\Phi$ & Design space \\ 
        $\mu$ & Mean value \\
        $\sigma$ & Elements of the observation error covariance matrix \\
        $\tilde{\sigma}$ & Elements of the inverse of the observation error covariance matrix \\
        $\boldsymbol{\Sigma}_y$ & Observation error covariance matrix \\
        $\boldsymbol{\Psi}$   & FIM metric \\
    \end{tabular}
\end{table}

\newpage
\bibliographystyle{elsarticle-harv} 
\bibliography{cas-refs}






\newpage

\appendix

\setcounter{equation}{0}
\renewcommand{\theequation}{S-\arabic{equation}}

\setcounter{figure}{0}
\renewcommand{\thefigure}{S-\arabic{figure}}

\setcounter{table}{0}
\renewcommand{\thetable}{S-\arabic{table}}

\setcounter{section}{0}
\renewcommand{\thesection}{S-\arabic{section}}

\setcounter{footnote}{0}

\setcounter{page}{1}
\renewcommand{\thepage}{S-\arabic{page}}

\newgeometry{left=1in, right=1in}

{\Large
\begin{center}
Supplementary Material for \\
 Measure This, Not That: \\
 Optimizing the Cost and Model-Based Information Content of Measurements
\end{center}
}

\begin{center}
Jialu Wang$^a$, Zedong Peng$^b$, Ryan Hughes$^{c,d}$, Debangsu Bhattacharyya$^{e}$, David E. Bernal Neira$^b$, 
Alexander W. Dowling\footnote{corresponding author: adowling@nd.edu}$^a$

\emph{$^a$Department of Chemical and Biomolecular Engineering, University of Notre Dame, Notre Dame, IN 46556}

\emph{$^b$Davidson School of Chemical Engineering, Purdue University, West Lafayette, IN 47907}

\emph{$^c$National Energy Technology Laboratory, Pittsburgh, PA 15236}

\emph{$^d$Support Contractor, National Energy Technology Laboratory, Pittsburgh, PA 15236}

\emph{$^e$Department of Chemical and Biomedical Engineering, West Virginia University, Morgantown, WV 26506}

\end{center}

\begin{center}
    June 12, 2024
\end{center}

\newpage

\section{Reaction Kinetics: Additional Data and Results}

\begin{table}[!htb]
\centering
\caption{Covariance structure for the kinetics case study.} 
\begin{tabular}{| l | cccccc |} \hline
& $C_A^{SCM}$ & $C_B^{SCM}$& $C_C^{SCM}$
& $C_A^{DCM}$ & $C_B^{DCM}$& $C_C^{DCM}$\\ \hline
$C_A^{SCM}$ & 1 & 0.1 & 0.1  & 0.5 & 0.05 & 0.05 \\
$C_B^{SCM}$ & 0.1 & 4 & 0.5  & 0.05 & 2 & 0.25 \\
$C_C^{SCM}$ & 0.1 & 0.5 & 8  & 0.05 & 0.25 & 4 \\
$C_A^{DCM}$ & 0.5 & 0.05 & 0.05 & 1 & 0.1 & 0.01 \\
$C_B^{DCM}$ & 0.05 & 2 & 0.25   & 0.1 & 4 & 0.5 \\
$C_C^{DCM}$ & 0.05 & 0.25 & 4   & 0.1 & 0.5 & 8 \\ \hline
\end{tabular}
\label{tab:kinetics_cov}
\end{table}

\begin{figure}[!ht]
\centering
\includegraphics[width=0.95\textwidth]{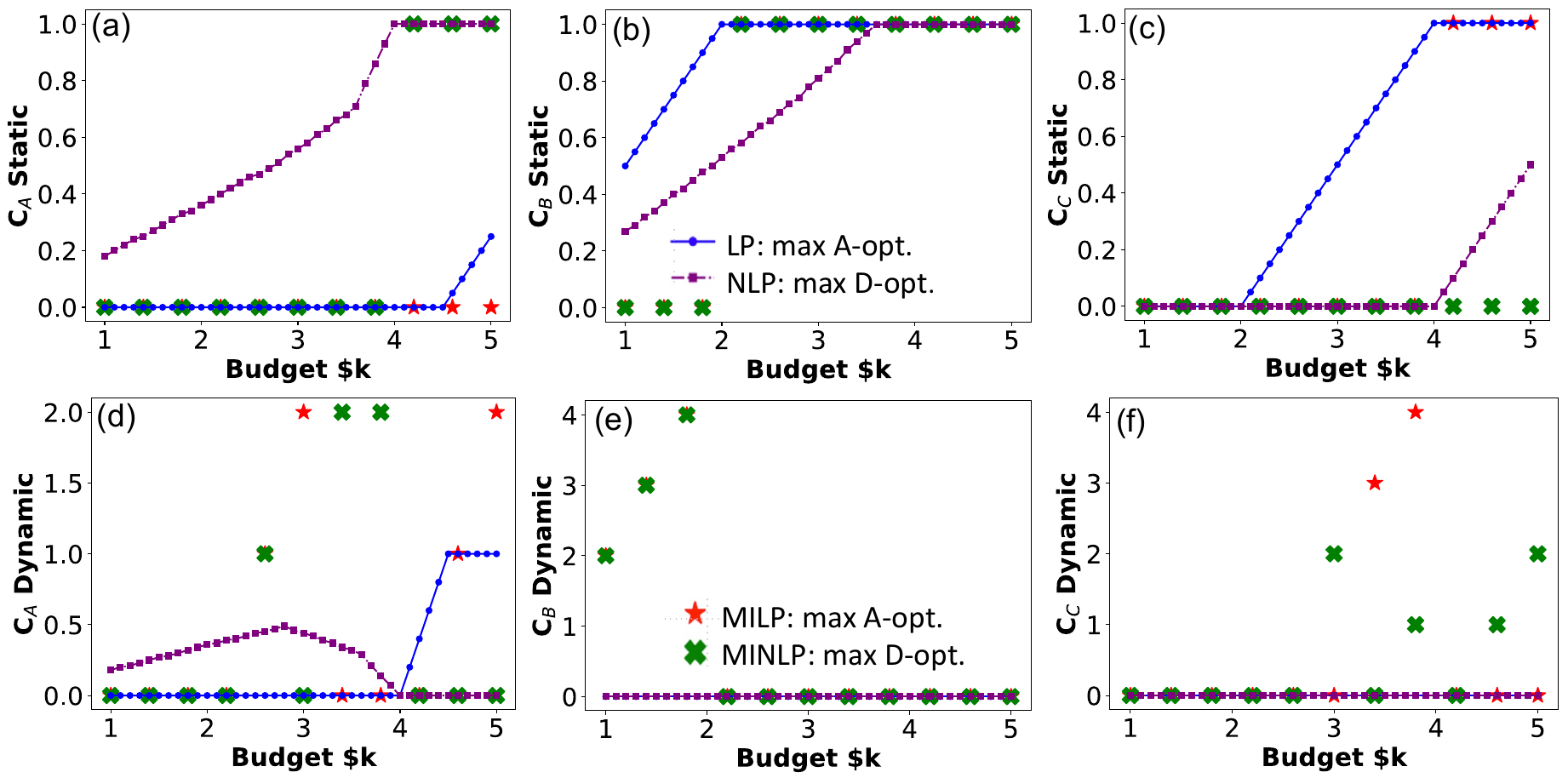}
\caption{Measurement solutions along with budgets in the kinetics case study. The top row contains (a), (b), (c), showing the binary decision 1 (yes) and 0 (no) of $C_A$, $C_B$, and $C_C$ as SCMs. The bottom row contains (d), (e), (f), showing the number of chosen time points of $C_A$, $C_B$, and $C_C$ as DCMs. 
Solutions of four optimization strategies are plotted: maximizing A-optimality of the MILP problem (red stars), its relaxed LP problem (blue line); maximizing D-optimality of the MINLP problem (green crosses), and its relaxed NLP problem (purple line). } 
\label{fig:kinetics_solution}
\end{figure}

\section{Reacton Kinetics: Consideration of the Operating Costs}
\label{sec:operate}

Operating and maintenance (O\&M) costs, also called direct operating costs, are an important part of evaluating investments and preparing budgets for chemical process industries \citep{garrett2012chemical}. The measurement optimization (MO) framework allows the inclusion of the operating costs of static-cost measurements (SCMs). We are showing the workflow to include the operating costs with the kinetics case study. 

A high-performance liquid chromatography (HPLC) with automated sample extraction processes can be integrated directly into the experimental system to continuously monitor the measurements as SCMs.
The instrument price is \$40,000, and its operating cost per year is approximately 15\% of the instrument cost, which is \$6,000 \citep{liptak2003instrument}. A reasonable assumption is the experimental team rents the access of an HPLC for this project for one year, and pays \$20 operating fee per use, which is equivalent to \$2.5 operating fee per time, which changes the cost structure to:

\begin{table}[htb!]
\centering
\caption{Measurement costs and types for the kinetics case study considering operating costs.}  
\begin{tabular}{lcc}
Name & Installation (\$) & Measure cost (\$/sample) \\
$C_A^{SCM}$ & 2000 & 2.5 \\
$C_B^{SCM}$ & 2000 & 2.5 \\
$C_C^{SCM}$ & 2000 & 2.5 \\
$C_A^{DCM}$ & 200 & 400 \\
$C_B^{DCM}$ & 200 & 400 \\
$C_C^{DCM}$ & 200 & 400 \\

\end{tabular}
\label{tab:kinetics_price_operate}
\end{table}

\begin{figure}[!ht]
\centering
\includegraphics[width=0.95\textwidth]{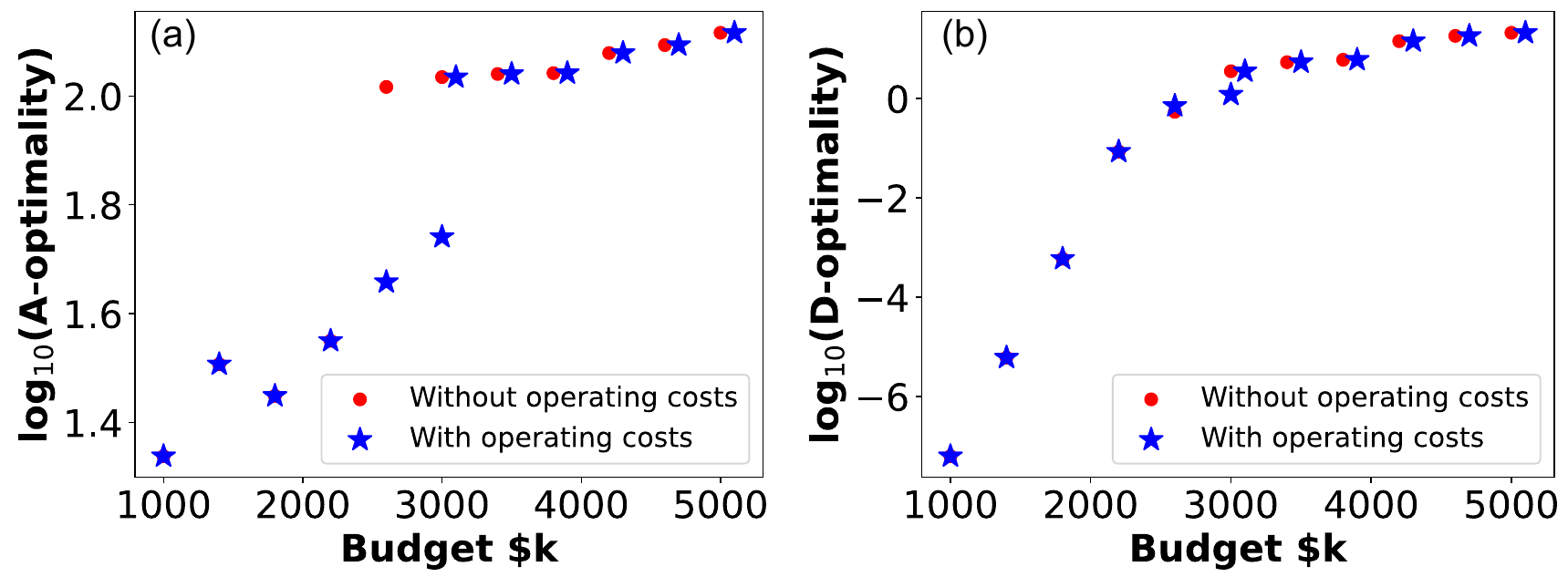}
\caption{Pareto-optimal trade-off between measurement budgets and A-optimality (trace of FIM) in (a), and the Pareto-optimal trade-off between measurement budgets and D-optimality (determinant of FIM) in (b) for the reaction kinetics case study when considering operating costs for SCMs. Results of maximizing D-optimality of the MINLP problem are plotted.} 
\label{fig:operate}
\end{figure}

Fig.~\ref{fig:operate} shows the comparison of the mixed-integer nonlinear programming (MINLP) problem maximizing D-optimality problem, with the cost strategy including operating costs, and the cost strategy excluding operating costs. When the budget is below \$2.6k, both strategies show the same results since there is no SCM chosen. When the budget is \$3k when not including operating costs, an SCM is chosen, while MO cannot afford this SCM when considering operating costs until the budget is increased to around \$3.1k. Generally, this problem chooses the same choices of measurements as the budget increases but needs \$20 more budget for each SCM it decides to choose. 

\FloatBarrier \vfill \pagebreak
\section{Reaction Kinetics: Consideration of Multi-component Samples}
\label{sec:multi}
Another practical scenario considered for the kinetics case study is that multi-component samples can allow the measurement of all three components, $C_A, C_B$, and $C_C$, with one sample. The cost strategy can be adapted to: 

\begin{table}[htb!]
\centering
\caption{Measurement costs and types for the kinetics case study considering operating costs.}  
\begin{tabular}{lcc}
Name & Installation (\$) & Measure cost (\$/sample) \\
$C^{SCM}$ & 2000 & 2.5 \\
$C^{DCM}$ & 200 & 400 \\
\end{tabular}
\label{tab:multi-component}
\end{table}

\begin{table}[htb!]
\centering
\caption{Optimal solutions for the D-optimality MINLPs for multi-component samples. For budgets of \$1k to \$5k, '1' means an SCM is selected and '0' means not selected. The selected time points [min] of DCMs are reported.}  
\begin{tabular}{|l|cccc|}
\toprule
Name & 1.0 & 1.4 & 1.8 & 2.2 \\
\midrule
$C_{A}^{SCM}$ & 0  & 0  & 0  & 1\\ 
$C_{B}^{SCM}$ & 0  & 0  & 0  & 1 \\
$C_{C}^{SCM}$ & 0  & 0  & 0  & 1 \\
\hline 
$C_{A}^{DCM}$ & 45  &  30  & 15  &  \\ 
& 60  &  45  & 30  &  \\  
&   &  60  & 45  &  \\  
&   &    & 60  &  \\  

$C_{B}^{DCM}$ & 45  &  30  & 15  &  \\ 
& 60  &  45  & 30  &  \\  
&   &  60  & 45  &  \\  
&   &    & 60  &  \\  

$C_{C}^{DCM}$ & 45  &  30  & 15  &  \\ 
& 60  &  45  & 30  &  \\  
&   &  60  & 45  &  \\  
&   &    & 60  &  \\  
\bottomrule
\end{tabular}
\label{tab:multi-component-result}
\end{table}

Table \ref{tab:multi-component-result} shows the optimization results maximizing D-optimality with the MINLP framework. With a budget of \$2.2k, all SCMs are chosen, and there is no need for more budgets. 

\FloatBarrier 
\vfill \pagebreak 

\section{Reaction Kinetics: Influence of Sampling Time Constraints}
\label{sec:sampling_time_constraint}

In many industrial context, the minimum time between samples and maximum number of samples within a working shift is an extremely important practical consideration. To show the flexibility of the proposed framework, we resolved the A-optimality MILP problem of the kinetics case study with a budget of \$3.4k. With our chosen time interval of 10 minutes in the paper, the selected DCMs include C$_C$ at 30, 45, and 60 minutes, resulting in A- and D-optimality values of 114.07 and 0.043, respectively. Increasing the time interval to 20 minutes, the selected DCMs include CC at 15, 37.5, and 60 minutes, and the A- and D-optimality values decrease to 112.69 and 0.096, respectively. The increase in the time interval (fewer possible sampling times) leads to a 1.22\% decrease in A-optimality, which is the objective function of the problem. The increase in D-optimality is because spreading out the DCMs results in selecting DCMs from a high-D-optimality region.

\section{Reaction Kinetics: Comparison of \texttt{CVXPY} and \texttt{Pyomo} results}
\label{sec:cvxpy}

\begin{figure}[!ht]
\centering
\includegraphics[width=0.95\textwidth]{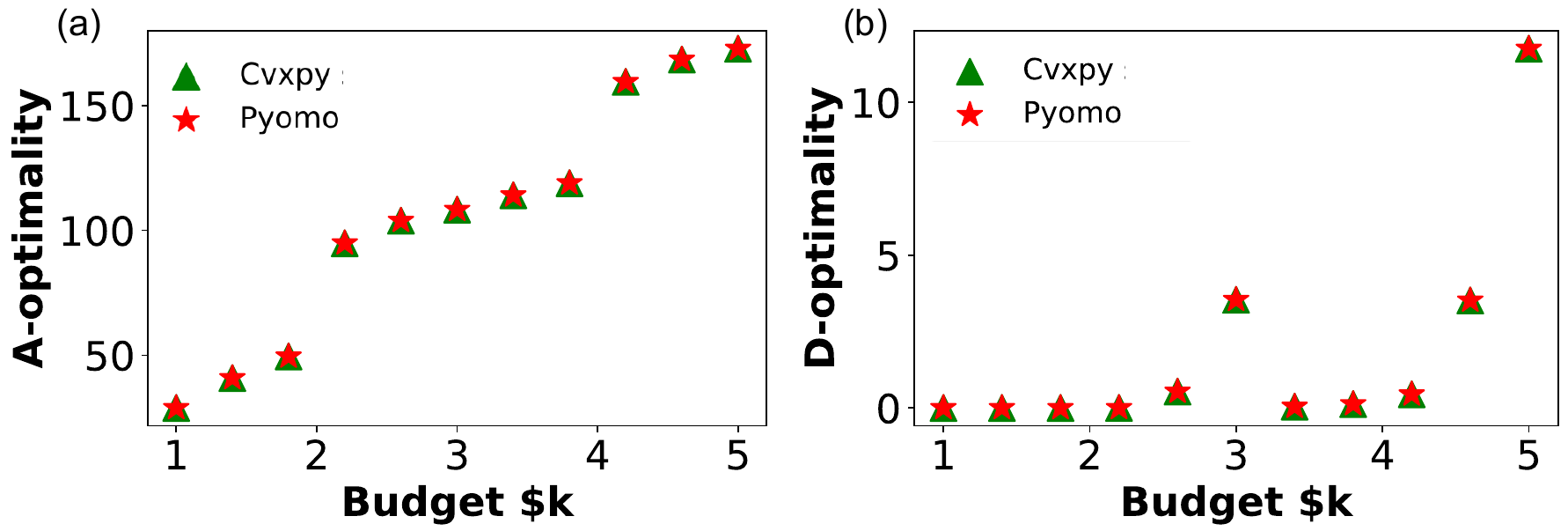}
\caption{Measurement solutions of the kinetics case study implemented in \texttt{CVXPY} and \texttt{Pyomo}, for the A-optimality MILP problem. } 
\label{fig:cvxpy_A}
\end{figure}

\begin{figure}[ht]
\centering
\includegraphics[width=0.95\textwidth]{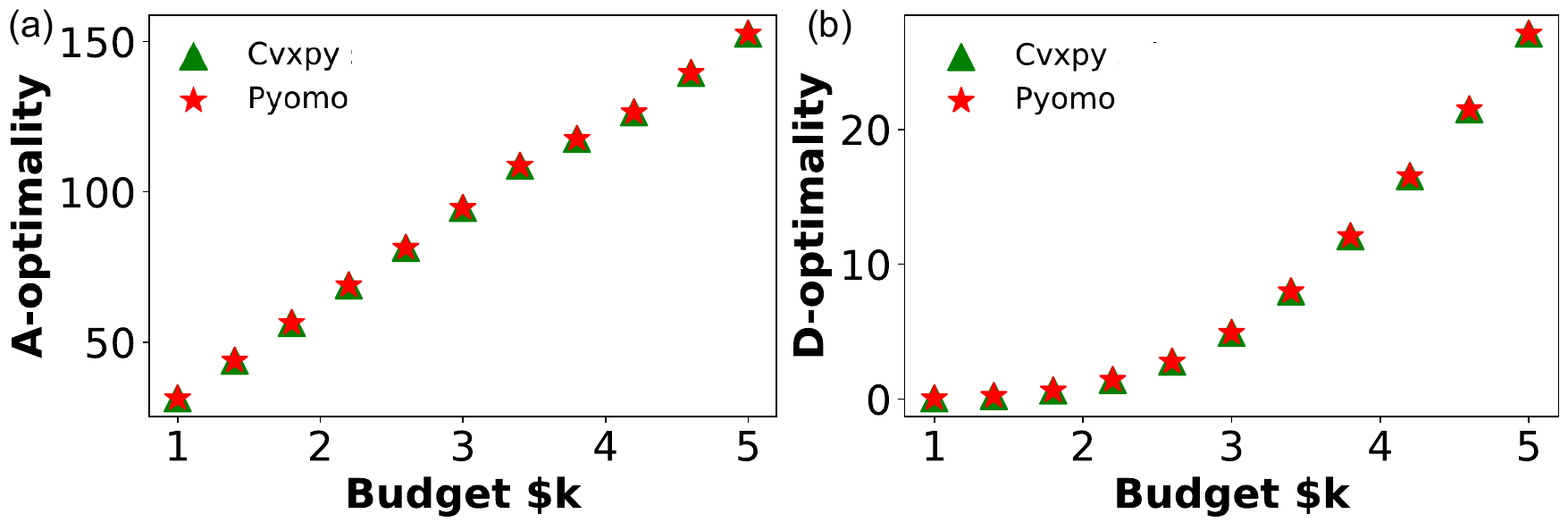}
\caption{Measurement solutions of the kinetics case study implemented in \texttt{CVXPY} and \texttt{Pyomo}, for the D-optimality NLP problem. } 
\label{fig:cvxpy_D_nlp}
\end{figure}

\FloatBarrier \vfill \pagebreak
\section{Rotary Bed: D-Optimality MINLP Solutions}
\label{sec:rotaryD}
\setlength\tabcolsep{3pt} 
\begin{longtable}{|l | c | c | c | c | c | c| c| c |c |c |c |c |c |c |c |c |c |c |c |c |c |c |c |c |c| }
\caption{Optimal solutions of SCMs for the D-optimality MINLPs where each column corresponds to a measurement budget of \$1l \$25k. '1' and '0' indicate if each SCM is selected or not selected. The name B. in the leftmost cell of the first row represents the row name, budget [\$k].}\\
\hline
\multicolumn{1}{|c|}{B.} & 
\multicolumn{1}{c|}{\textbf{1}} & \multicolumn{1}{c|}{\textbf{2}} & \multicolumn{1}{c|}{\textbf{3}} &
\multicolumn{1}{c|}{\textbf{4}} &
\multicolumn{1}{c|}{\textbf{5}} &
\multicolumn{1}{c|}{\textbf{6}} &
\multicolumn{1}{c|}{\textbf{7}} &
\multicolumn{1}{c|}{\textbf{8}} &
\multicolumn{1}{c|}{\textbf{9}} &
\multicolumn{1}{c|}{\textbf{10}} &
\multicolumn{1}{c|}{\textbf{11}} &
\multicolumn{1}{c|}{\textbf{12}}&
\multicolumn{1}{c|}{\textbf{13}}&
\multicolumn{1}{c|}{\textbf{14}}&
\multicolumn{1}{c|}{\textbf{15}}&
\multicolumn{1}{c|}{\textbf{16}}&
\multicolumn{1}{c|}{\textbf{17}}&
\multicolumn{1}{c|}{\textbf{18}}&
\multicolumn{1}{c|}{\textbf{19}}&
\multicolumn{1}{c|}{\textbf{20}}&
\multicolumn{1}{c|}{\textbf{21}}&
\multicolumn{1}{c|}{\textbf{22}}&
\multicolumn{1}{c|}{\textbf{23}}&
\multicolumn{1}{c|}{\textbf{24}}&
\multicolumn{1}{c|}{\textbf{25}}\\ \hline 
\endfirsthead

\multicolumn{26}{c}%
{{\bfseries \tablename\ \thetable{} -- continued from previous page}} \\
\hline 
\multicolumn{1}{|c|}{B.} & 
\multicolumn{1}{c|}{\textbf{1}} & \multicolumn{1}{c|}{\textbf{2}} & \multicolumn{1}{c|}{\textbf{3}} &
\multicolumn{1}{c|}{\textbf{4}} &
\multicolumn{1}{c|}{\textbf{5}} &
\multicolumn{1}{c|}{\textbf{6}} &
\multicolumn{1}{c|}{\textbf{7}} &
\multicolumn{1}{c|}{\textbf{8}} &
\multicolumn{1}{c|}{\textbf{9}} &
\multicolumn{1}{c|}{\textbf{10}} &
\multicolumn{1}{c|}{\textbf{11}} &
\multicolumn{1}{c|}{\textbf{12}}&
\multicolumn{1}{c|}{\textbf{13}}&
\multicolumn{1}{c|}{\textbf{14}}&
\multicolumn{1}{c|}{\textbf{15}}&
\multicolumn{1}{c|}{\textbf{16}}&
\multicolumn{1}{c|}{\textbf{17}}&
\multicolumn{1}{c|}{\textbf{18}}&
\multicolumn{1}{c|}{\textbf{19}}&
\multicolumn{1}{c|}{\textbf{20}}&
\multicolumn{1}{c|}{\textbf{21}}&
\multicolumn{1}{c|}{\textbf{22}}&
\multicolumn{1}{c|}{\textbf{23}}&
\multicolumn{1}{c|}{\textbf{24}}&
\multicolumn{1}{c|}{\textbf{25}}\\ \hline 
\endhead

\hline \multicolumn{26}{|r|}{{Continued on next page}} \\ \hline
\endfoot

\hline \hline
\endlastfoot
 
$F_{in}^{ads}$ & 0  & 0  & 0  & 0 & 1
& 1& 1 & 1& 1& 1
& 1& 1& 1& 1& 1
& 1& 1& 1& 1& 1
& 1 & 1& 1& 1& 1 \\ 
$F_{out}^{ads}$ & 0  & 0  & 0  & 0 & 0
& 0& 0 & 0& 0& 1
& 1& 0& 0& 0& 0
& 0& 1& 1& 1& 1
& 0 & 0& 1& 1& 1 \\
$T_{out}^{ads}$ & 1  & 1  & 1  & 1 & 1
& 1& 1 & 1& 1& 1
& 1& 1& 1& 1& 1
& 1& 1& 1& 1& 1
& 1 & 1& 1& 1& 1 \\
$F_{in}^{des}$ & 0  & 0  & 0  & 0 & 0
& 1& 1 & 1& 1 & 1
& 1& 1& 1& 1 & 1
& 1& 1& 1& 1 & 1
& 1 & 1& 1& 1& 1 \\
$F_{out}^{des}$ & 0  &1  & 1  & 1 & 1& 1& 1
& 1& 1& 1& 1& 1& 1& 1& 1& 1& 1& 1& 1& 1& 1 
& 1& 1& 1& 1 \\
$T_{out}^{des}$ & 1  & 1  & 1  & 1 & 1& 1& 1
& 1& 1& 1& 1& 1& 1& 1& 1& 1& 1& 1& 1& 1& 1 
& 1& 1& 1& 1 \\
$T_{19}^{ads}$ & 0  & 0  & 1  & 1 & 1& 1& 1
& 1& 1& 1& 1& 1& 1& 1& 1& 1& 1& 1& 1& 1& 1 
& 1& 1& 1& 1 \\
$T_{23}^{ads}$ & 0  & 0  & 0  & 1 & 1& 1& 1
& 1& 1& 1& 1& 0& 1& 1& 1& 1& 1& 1& 1& 1& 1 
& 1& 1& 1& 1 \\
$T_{28}^{ads}$ & 0  & 0  & 0  & 0 & 0& 0& 0
& 1& 1& 1& 1& 0& 0& 1& 1& 1& 1& 1& 1& 0
& 0 & 1& 1& 1& 1 \\
$z_{out}^{ads}$ & 0  & 0  & 0  & 0 & 0& 0& 0
& 0& 0& 0& 0& 1& 1& 1& 1& 1& 1& 1& 1& 1
& 1  & 1& 1& 1& 1 \\
$z_{out}^{des}$ & 0  & 0  & 0  & 0 & 0& 0& 0
& 0& 0& 0& 0& 0& 0& 0& 0& 0& 0& 0& 0& 1
& 1 & 1& 1& 1& 1 \\
\end{longtable}
\label{tab:rotary_solution_SCM_D}

\setlength\tabcolsep{3pt} 
\begin{longtable}{|l | c | c | c | c | c | c| c| c |c |c |c |c |c |c |c |c |c |c |c |c |c |c |c |c |c| }
\caption{Optimal solutions of DCMs for the D-optimality MINLPs where each column corresponds to a measurement budget of \$1k to \$25k. The selected sample times [min] of DCMs are reported. The name B. in the leftmost cell of the first row represents the row name, budget [\$k].}\\
\hline
\multicolumn{1}{|c|}{B.} & 
\multicolumn{1}{c|}{\textbf{7}} &
\multicolumn{1}{c|}{\textbf{8}} &
\multicolumn{1}{c|}{\textbf{9}} &
\multicolumn{1}{c|}{\textbf{10}} &
\multicolumn{1}{c|}{\textbf{11}} &
\multicolumn{1}{c|}{\textbf{12}} &
\multicolumn{1}{c|}{\textbf{14}}&
\multicolumn{1}{c|}{\textbf{15}}&
\multicolumn{1}{c|}{\textbf{16}}&
\multicolumn{1}{c|}{\textbf{17}}&
\multicolumn{1}{c|}{\textbf{18}}&
\multicolumn{1}{c|}{\textbf{19}}&
\multicolumn{1}{c|}{\textbf{20}}&
\multicolumn{1}{c|}{\textbf{21}}&
\multicolumn{1}{c|}{\textbf{22}}&
\multicolumn{1}{c|}{\textbf{23}}&
\multicolumn{1}{c|}{\textbf{24}}&
\multicolumn{1}{c|}{\textbf{25}}\\ \hline 
\endfirsthead

\multicolumn{26}{c}%
{{\bfseries \tablename\ \thetable{} -- continued from previous page}} \\
\hline 
\multicolumn{1}{|c|}{B.} & 
\multicolumn{1}{c|}{\textbf{7}} &
\multicolumn{1}{c|}{\textbf{8}} &
\multicolumn{1}{c|}{\textbf{9}} &
\multicolumn{1}{c|}{\textbf{10}} &
\multicolumn{1}{c|}{\textbf{11}} &
\multicolumn{1}{c|}{\textbf{12}} &
\multicolumn{1}{c|}{\textbf{14}}&
\multicolumn{1}{c|}{\textbf{15}}&
\multicolumn{1}{c|}{\textbf{16}}&
\multicolumn{1}{c|}{\textbf{17}}&
\multicolumn{1}{c|}{\textbf{18}}&
\multicolumn{1}{c|}{\textbf{19}}&
\multicolumn{1}{c|}{\textbf{20}}&
\multicolumn{1}{c|}{\textbf{21}}&
\multicolumn{1}{c|}{\textbf{22}}&
\multicolumn{1}{c|}{\textbf{23}}&
\multicolumn{1}{c|}{\textbf{24}}&
\multicolumn{1}{c|}{\textbf{25}}\\ \hline 
\endhead

\hline \multicolumn{26}{|r|}{{Continued on next page}} \\ \hline
\endfoot

\hline \hline
\endlastfoot
 
$z_{out}^{ads}$ &\ &\ &126 &126 
&124  &36  &\  &\ 
&\ &\  &\ &\ 
&\ &\ &\ &\ &\  &\ 
\\
 &\ &\ &136 &136 
 &134 &126 &\ &\ 
 &\ &\ &\ &\  
 &\ &\ &\ &\ &\  &\ 
\\
 &\ &\ &180 &182 
 &180 &136  &\  &\ 
 &\ &\ &\ &\  
 &\ &\ &\ &\ &\  &\ 
\\ &\ &\ &190 &192 
&190 & 186  &\  &\ 
&\ &\ &\ &\  
&\ &\ &\ &\ &\  &\ 
\\
&\ &\ &\ &\ 
&200 &196 &\ &\ 
&\ &\ &\ &\  
&\ &\ &\ &\ &\  &\ 
\\
$z_{out}^{des}$   &\ 
 &\ &2 &2 &2 &2  &\ 
 &\ &2 &2 &2 &12  &12 &\
 &\ &\ &\ &\ \\

 &\ 
 &\ &12 &12 &12 &12  &\ 
   &\ &12 &12 &12 &166 &166  &\ 
   &\ &\ &\ &\ \\

 &\ 
 &\ &22 &22 &22 &24  &\
  &\ &22 &22 &22 &176 &176 &\
  &\ &\ &\ &\
  \\

&\ 
 &\ &32 &32 &32  &156  &\ 
  &\ &32 &32 &32 &186  &186 &\
  &\ &\ &\ &\
   \\

 &\ 
 &\ &42 &42 &42 &166 &\ 
  &\ &42 &42 &42 &196 &196  &\
  &\ &\ &\ &\
  \\

$z_{19}^{ads}$   
&\ &\ &\ &\ 
&\ &146 &\ &\ 
&\ &\ &158 &2  &2
&\   &\ &\ & 168 &76 \\

 &\ &\ &\ &\ 
 &\ &\ &\ &\ 
 &\ &\ &186  &50  &50
 &\  &\ &\ & 180 &140 \\

 &\ &\ &\ &\ 
 &\ &\  &\ &\ 
 &\ &\ &196  &136 &136
 &\ &\ &\ & 196 &150\\

 &\ &\ &\ &\ 
 &\ &\ &\ &\ 
 &\ &\ &206  &146  &146
 &\  &\ &\ & 206 &180\\

 &\ &\ &\ &\ 
 &\ &\ &\ &\ 
 &\ &\ &216  &156  &156
 &\ &\ &\ & 216 &190\\

$z_{23}^{ads}$   &\ 
 &\ &\ &\ &158 &66   &\ 
  &\ &\ &\ &\ &22 &22 &\ 
   &\ &\ &\ &160\\

&\ 
  &\ &\ &\ &210 &76 &\ 
   &\ &\ &\ &\ &66 &66 &\ 
    &\ &\ &\ &170 \\

  &\ 
   &\ &\ &\ &220  &206  &\ 
    &\ &\ &\ &\ &76 &76  &\ 
     &\ &\ &\ &200 \\

&\ 
    &\ &\ &\ &\ &216  &\ 
     &\ &\ &\ &\ &206 &206 &\ 
      &\ &\ &\ &210 \\

&\ 
     &\ &\ &\ &\  &\ 
      &\ &\ &\ &\    &\ &216 &216  &\ 
    &\ &\ &\ &220 \\
    
$z_{28}^{ads}$  &86 
&86 &86 &86 &84 &86 &86 
&86 &86 &86 &86 &86 &86  &86
&86 &86 &86 &86\\

 &96 
 &96 &96 &96 &94  &96 &96
 &96 &96 &96 &96 &96 &96 &96
 &96 &96 &96 &96\\

 &106 
  &106 &106 &106 &104  &106 &106
  &106 &106 &106 &106 &106 &106  &106
   &106 &106 &106 &106\\

  &116 
   &116 &116 &116 &114 &116 &116
   &116 &116 &116 &116 &116 &116 &116
   &116 &116 &116 &116\\

   &126 
    &126 &\ &\ &170 &176 &126 
    &126 &126 &126 &126 &126 &126 &126
    &126 &126 &126 &126  \\
\end{longtable}
\label{tab:rotary_solution_DCM_D}

\FloatBarrier \vfill \pagebreak

\section{Rotary Bed: Further Analysis of A-Optimality Results}
\label{sec:rotaryA}

In the main text, Fig. \ref{fig:rotary}, Table \ref{tab:rotary_solution_SCM_A}, and Table \ref{tab:rotary_solution_DCM_A} show the best A-optimality results for maximum experimental budgets ranging from \$1k to \$25k for the rotary bed case study. We now further analyze these results. Careful inspection shows that the proposed optimization framework selects the measurements in order from most to least informative. With a budget of less than \$6k, only static-cost measurements are chosen as a balance of budget and information content. Starting with \$1k, MO chooses two temperature static measurements, adsorption and desorption gas outlet temperature, with a cost of \$0.5k each. Starting with \$2k, every time \$1k is added to the budget, a \$1k static measurement is able to be added. At \$2k, desorption gas outlet flowrate is added. At \$3k, desorption inlet flowrate is added. At \$4k, adsorption gas inlet flowrate is added. At \$5k, the adsorption temperature inside bed at position 19 is added. At \$6k, the adsorption temperature inside the bed at position 23 is added. 

These optimal solutions can be verified by the trace of single measurements. A-optimality, as the trace of the Fisher information matrix (FIM), is the sum of diagonal elements of FIM. The sum of A-optimality is, therefore, a linear equation of the A-optimality of unit FIMs. Without correlations between measurements, A-optimality is only computed by the unit FIMs representing the variance of one measurement. Desorption gas outlet temperature is firstly chosen at a budget of \$1k because it has the biggest A-optimality value of 31328. Adsorption gas outlet temperature is also chosen, although it has a lower A-optimality of 154, since the budget only allows another \$0.5k measurement, and it is the only remaining one. Starting from \$2k, all remaining static-cost measurements cost \$1k; therefore, they are chosen according only to their variance FIM. Desorption gas outlet flowrate has an A-optimality of 627 and is chosen for \$2k. The next big A-optimal measurement is desorption gas inlet flowrate, which is chosen for \$3k. In this way, the optimal solution for the first \$6k budget can be verified. 

Starting from \$7k, MO begins to add dynamic-cost measurements. With the constraint of choosing, at most, 5 time points for each dynamic-cost measurement, MO chooses 5 time points from the adsorption mole fraction at position 28 with a \$7k budget. When the budget increases to \$11 k, MO chooses 18 time points from all dynamic-cost measurements in total, under the constraint that the total time points should be under 20, and the time points should be at least 10 minutes from each other. 

For budgets of \$11k and \$12k, MO removes the selection of 3 static-cost measurements and all time points of dynamic-cost measurements to allow the budget for selecting an expensive static-cost measurement, gas chromatography (GC) in the adsorption gas outlet, which costs \$7k. With the budget increasing to \$19k, MO, while persistently selecting the GC in the adsorption gas outlet end, adds back the static-cost measurements and time points of dynamic-cost measurements step by step, constrained by the budget. At \$19k, all static-cost measurements are selected except for the other expensive GC of the desorption gas outlet, and 20 time points in total for dynamic-cost measurements are selected. Since GC in the adsorption gas outlet is selected as a static-cost measurement now, no time points for its dynamic-cost option are chosen; instead, 5 time points each are chosen for the other four dynamic-cost measurements. 

For budgets of \$19 to \$20k, MO chooses to remove the selection of all dynamic-cost measurements and some of the static-cost measurements, to include the other expensive static-cost choice, GC of desorption gas outlet. For a \$20k to \$25k budget, MO adds back these static-cost measurements and time points for dynamic-cost measurements step by step. At a \$25k budget, all static-cost measurements are selected. Mole fractions of adsorption and desorption gas outlet are chosen as static-cost measurements, hence no time points for its dynamic-cost options are chosen. 5 time points each are chosen from all the remaining three dynamic-cost measurements. 

In summary, for this case study, MO tends to prioritize static-cost measurements firstly with a small budget, then fills the budget with dynamic measurements. It can be observed that the temperature of the adsorption and desorption gas outlet end are chosen for their high information content contributions and low costs for every budget. On the contrary, the flowrate of the adsorption gas outlet end is considered by MO as a comparatively uninformative measurement, which is only chosen when there is not enough budget to choose other measurements. The two GC of adsorption and desorption, although informative, also cost as high as static-cost measurements. When not having enough budget, MO selects them as dynamic-cost measurements to provide information while costing low. Starting from \$12k, MO manages to include at least one of the GC machines as a static-cost measurement at the cost of removing several other measurements.

\end{document}